\documentclass[reqno]{amsart}
\numberwithin{equation}{section}
\usepackage{graphicx}

\newtheorem{theorem}{\sc Theorem}[section]
\newtheorem{lemma}{\sc Lemma}[section]

\newtheorem{remark}{\sc Remark}

\def\bdy #1{\partial #1}

\def\bbR{{\mathbb R}}

\def\bbE{{\mathcal E}}

\def\div{\operatorname{div}}
\def\Def{\operatorname{Def}}
\def\id{{\text{Id}}}
\def\supp{{\text{supp}}}
\def\L{{\mathcal L}}
\def\ta{{\tilde{a}}}

\def\tb{{\tilde{b}}}

\def\tv{{\tilde{v}}}

\def\ttv{{\bar{v}}}
\def\tq{{\tilde{q}}}
\def\tg{{\tilde{g}}}

\def\tn{{\tilde{n}}}

\def\tL{{\tilde{\mathcal L}}}
\def\tF{{\tilde{F}}}
\def\teta{{\tilde{\eta}}}

\def\nv{{v_{\kappa}}}
\def\na{{a_{\kappa}}}

\def\nq{{q_{\kappa}}}
\def\neta{{\eta_{\kappa}}}
\def\nn{{n_{\kappa}}}
\def\nb{{b_{\kappa}}}
\def\ng{{g_{\kappa}}}
\def\nX{{X_{\kappa}}}
\def\mP{{{\mathcal P}(\|\eta_\kappa\|^2_{H^{2.5}(\Gamma)})}}

\def\H1H2{{H^{1;2}(\Omega;\Gamma)}}

\def\aa{\sqrt{a}a^{\alpha\beta\gamma\delta}}

\begin{document}
\title[Navier-Stokes interacting with a solid shell]
{Navier-Stokes equations interacting with a nonlinear elastic
solid shell}

\author{C.H. Arthur Cheng}
\email{cchsiao@math.ucdavis.edu}
\author{Daniel Coutand}
\email{coutand@math.ucdavis.edu}
\author{Steve Shkoller}
\email{shkoller@math.ucdavis.edu}

\address{Department of Mathematics, University of California, Davis, CA 95616}

\maketitle
\begin{abstract}
We study a moving boundary value problem consisting of a viscous
incompressible fluid moving and interacting with a nonlinear elastic
solid shell.  The fluid motion is governed by the Navier-Stokes equations,
while the shell is modeled by the nonlinear Koiter shell model, consisting
of both bending and membrane tractions.
The fluid is coupled to the solid shell through
continuity of displacements and tractions (stresses)
along the moving material interface.  We prove existence and
uniqueness of solutions in Sobolev spaces.
\end{abstract}

\section{Introduction}
\subsection{The problem statement and background}
Fluid-solid interaction problems involving moving material
interfaces have been the focus of active research since the
nineties. The first problem solved in this area was for the case of
a rigid body moving in a viscous fluid (see \cite{DeEs}, \cite{GrMa}
and also the early works of \cite{Wein} and \cite{Serre} for a rigid
body moving in a Stokes flow in the full space). The case of an
elastic body moving in a viscous fluid was considerably more
challenging because of some apparent regularity incompatibilities
between the parabolic fluid phase and the hyperbolic solid phase.
The first existence results in this area were for regularized
elasticity laws, such as in \cite{DeEsGrTa} for a {\it finite}
number of elastic modes, or in \cite{DaVe}, \cite{ChDeEsGr}, and
\cite{Bo2005} for hyperviscous elasticity laws, or in \cite{LiWa} in
which a phase-field regularization ``fattens'' the sharp interface
via a diffuse-interface model.

The treatment of classical elasticity laws for the solid phase,
without any regularizing term, was only considered recently in
\cite{CoSh2005} for the three-dimensional linear St.
Venant-Kirchhoff constitutive law and in \cite{CoSh2006} for
quasilinear elastodynamics coupled to the Navier-Stokes equations.
Some of the basic new ideas introduced in those works concerned a
functional framework that scales in a hyperbolic fashion (and is
therefore driven by the solid phase), the introduction of
approximate problems either penalized with respect to the
divergence-free constraint in the moving fluid domain, or smoothed
by an appropriate parabolic artificial viscosity in the solid phase
(chosen in an asymptotically convergent and consistent fashion), and
the tracking of the motion of the interface by difference quotients
techniques.

The complimentary fluid-solid interaction problem, studied herein,
consists of the motion of
a viscous incompressible fluid enclosed by
a moving thin nonlinear elastic {\it solid} shell.   Our companion
paper \cite{ChCoSh2006} treats the case of a viscous incompressible fluid
enclosed by a moving thin nonlinear elastic {\it fluid} shell.
This is a moving
boundary problem that models the motion of a viscous incompressible
Newtonian fluid inside of a deformable elastic structure. The main
mathematical differences with respect to the previous problem of a
deformable solid body moving inside of the fluid is that the
shell encloses the fluid and is mathematically the {\it boundary} of the
fluid.  The shell model consists of ``elliptic'' operators which do not
provide the expected regurality associated with the highest order operator
coming from the shell's bending energy, and, in particular,  ellipticity holds
only for short time.  The only cases considered until now consisted
of regularized problems, wherein the elliptic degeneracy occurs along
a {\it fixed direction}, such as in \cite{FlOr} or \cite{ChDeEsGr}.

We are concerned here with establishing the existence and uniqueness
of solutions to the time-dependent incompressible Navier-Stokes
equations interacting with a quasilinear elastic solid shell of Koiter
type (see \cite{PGC} for a detailed account of Koiter shells).
The solid shell energy is a nonlinear function of the
first and second fundamental forms of the moving boundary.

Let $\Omega\subset\bbR^n$ denote an open bounded domain with
boundary $\Gamma := \bdy\Omega$. For each $t\in (0,T]$, we wish to
find the domain $\Omega(t)$, a divergence-free velocity field
$u(t,\cdot)$, a pressure function $p(t,\cdot)$ on $\Omega(t)$, and a
volume-preserving transformation $\eta(t,\cdot):\Omega \to \bbR^n$
such that
\begin{subequations}\label{NSequation}
\begin{alignat}{2}
\Omega(t) &= \eta(t,\Omega) \,, & \\
\eta_t(t,x) &= u(t,\eta(t,x)) \,, & \\
u_t + \nabla_u u - \nu\Delta u &= -\nabla p + f & \text{in }\
\Omega(t)\,,
\label{NSequation.a} \\
\div u &= 0 & \text{in }\ \Omega(t)\,,
\label{NSequation.b}\\
(\nu\Def u -p\id)n &= {\mathfrak t}_{shell}  & \ \ \ \text{on } \
\Gamma(t)\,,
\label{NSequation.c} \\
u(0,x) &= u_0(x) &\forall x\in \Omega \,, \\
\eta(0,x) &= x &\forall x\in \Omega \,,
\end{alignat}
\end{subequations}
where $\nu$ is the kinematic viscosity, $n(t,\cdot)$ is the outward
pointing unit normal to $\Gamma(t)$, $\Gamma(t) := \partial
\Omega(t)$ denotes the boundary of $\Omega(t)$, $\Def u$ is twice
the rate of deformation tensor of $u$, given in coordinates by
$u^i_{,j} + u^j_{,i}$, where $u^i_{,j}$ denotes $\frac{\partial u^i}
{\partial x^j}$, and ${\mathfrak t}_{shell}$ is the traction
imparted onto the fluid by the elastic solid shell, which we describe
next.

With $\varepsilon$ denoting the thickness of the Koiter shell, $\lambda/2$ and
$\mu/2$ the Lam$\acute{\text{e}}$ constants, the energy stored in
the elastic surface has the form
\begin{align*}
\varepsilon E_{mem} + \frac{\varepsilon^3}{3} E_{ben}\,,
\end{align*}
where the membrane energy $E_{mem}$ is
\begin{align}
E_{mem} = \int_{\Gamma} \aa (g_{\alpha\beta} - g_{0\alpha\beta})
(g_{\gamma\delta} - g_{0\gamma\delta}) dS_0 \label{membrane}
\end{align}
and the bending energy $E_{ben}$ is
\begin{align}
E_{ben} = \int_{\Gamma} \aa (b_{\alpha\beta} - b_{0\alpha\beta})
(b_{\gamma\delta} - b_{0\gamma\delta}) dS_0\,, \label{bending}
\end{align}
where
\begin{align*}
a^{\alpha\beta\gamma\delta} = \frac{4\lambda\mu}{\lambda + 2\mu}
g_0^{\alpha\beta} g_0^{\gamma\delta} + 2\mu (g_0^{\alpha\gamma}
g_0^{\beta\delta} + g_0^{\alpha\delta} g_0^{\beta\gamma})\,.
\end{align*}
We let
$$g_{\alpha\beta} = \eta_{,\alpha}\cdot \eta_{,\beta}
 \ \ \text{ denote the induced metric on } \ \Gamma(t),$$
and
$$b_{\alpha\beta} = \eta_{,\alpha\beta}\cdot n
 \ \ \text{ denotes the second fundamental form.} $$
$g_0$ and $b_0$ denote the induced metric and
second fundamental form of the unstressed initial configuration at $t=0$.

The traction vector
$${\mathfrak t}_{shell} = \varepsilon {\mathfrak t}_{mem} + \frac{\varepsilon^3}{3} {\mathfrak t}_{ben}$$
is computed from the first variation of the energy function $E_{shell}$, and
will be stated in Section \ref{explicitform}.

In this paper, we will prove well-posedness for this system in
the case that the fluid is two-dimensional and the solid shell is its
one-dimensional closed boundary.

\section{Formulation of the problem}\label{explicitform}
\subsection{Fundamental geometric identities}
We use $(\cdot)'$ to denote the derivative of $(\cdot)$ along the boundary, and we
use $\delta \eta$ to denote the first variation of $\eta$.
The following formulas will be used often:
\begin{subequations}\label{normalderivatives}
\begin{align}
\delta n &= - |\eta'|^{-2} (n\cdot\delta\eta')\eta' \,,\\
n' &= -|\eta'|^{-2} (\eta''\cdot n)\eta'= -g^{-1/2} b \tau \,,\\
n'' &= 3 |\eta'|^{-4}(\eta''\cdot\eta')(\eta''\cdot n)\eta' -
|\eta'|^{-2}(\eta'''\cdot n)\eta' - |\eta'|^{-2} (\eta''\cdot n)\eta
'' \nonumber\\
&= - |\eta'|^{-2} (\eta''\cdot n)^2 n + |\eta'|^{-4}
(\eta''\cdot\eta') (\eta''\cdot n)\eta' - |\eta'|^{-2} (\eta''\cdot
n)' \eta' \nonumber\\
&= -g^{-1}b^2 n + \Big[\frac{g^{-3/2}}{2} g' b - g^{-1/2} b'\Big] \tau \,, \\
\eta'' &= (\eta''\cdot n)n + |\eta'|^{-1}(\eta''\cdot\eta')\tau = bn
+ \frac{g'}{2\sqrt{g}} \tau\,,
\end{align}
\end{subequations}
where $g = |\eta'|^2$, $\tau = g^{-1/2}\eta'$, and $b = \eta''\cdot
n$.

\subsection{The shell traction}

The bending energy (\ref{bending}) and membrane energy (\ref{membrane}) are
expressed as
$$
E_{ben} = \int_\Gamma |\eta_0'|^{-3} (b-b_0)^2 dS_0 \,, \ \ \
E_{mem} = \int_\Gamma  |\eta_0'|^{-3} (g-g_0)^2 dS_0 \,.
$$
Computing the first variation of the bending energy, we find that
the bending traction $\L_b$ is given by
\begin{align*}
\L_b(\eta) = 2 \Big[|\eta_0'|^{-3} (b - b_0)n\Big]'' +
\Big[|\eta_0'|^{-3} g^{-1}g' (b - b_0) n\Big]'\,,
\end{align*}
where, once again, $b_0 = \eta_0''\cdot n_0$ is the second fundamental form of
the unstressed initial boundary.

Taking the variation of $E_{mem}$ we find the
membrane traction $\L_m$  is
\begin{align*}
\L_m(\eta) = - 4 \Big[|\eta_0'|^{-3}(|\eta'|^2 -
|\eta_0'|^2)\eta'\Big]'\,.
\end{align*}

\subsection{Lagrangian formulation}  Let $\eta(t,x)=x+\int_0^t
u(s,x)ds$ denote the Lagrangian particle placement field, a
volume-preserving embedding of $\Omega$ onto $\Omega(t)\subset
{\mathbb R}^2$,  and denote the inverse matrix of $\nabla
\eta(x,t)$ by
\begin{align}
a(x,t) = [\nabla\eta(x,t)]^{-1} \label{adefn} \,.
\end{align}
Let $v=u\circ\eta$ denote the Lagrangian or material velocity field,
$q=p\circ\eta$ the Lagrangian pressure function, and $F=f\circ\eta$
the forcing function in the material frame.
The coupled fluid-structure problem has the following Lagrangian description:
\begin{subequations}\label{solidshelleq}
\begin{alignat}{2}
v=&\ \eta_t && \text{in }(0,T)\times\Omega\,,\\
v^i_t - \nu (a_\ell^j a_\ell^k v_{,k}^i)_{,j} =& - (a_i^j q)_{,j} + F^i && \text{in }(0,T)\times\Omega\,,\\
a_i^j v_{,j}^i =&\ 0 && \text{in }(0,T)\times\Omega \,,\\
\Big[\nu (a_i^k v_{,k}^j + a_j^k v^i_{,k})-q\delta_{ij}\Big]a_j^\ell
N_\ell =&\ \varepsilon \L_b(\eta) + \frac{\varepsilon^3}{3}
\L_m(\eta) \qquad&& \text{on }(0,T)\times\Gamma \,,\\
v(0,x)=&\ u_0(x) && \text{on } \{t=0\}\times\Omega\,,\\
\eta =&\ \id && \text{on } \{t=0\}\times\Omega\,.
\end{alignat}
\end{subequations}

\section{Notation and conventions}\label{notation}

\noindent For $T>0$, we set
\begin{align*}
{\mathcal V}^1(T) =&\ \Big\{ v\in L^2(0,T;H^1(\Omega))\ \Big|\ v_t \in L^2(0,T;H^1(\Omega)')\Big\}\,; \\
{\mathcal V}^2(T) =&\ \Big\{ v\in L^2(0,T;H^2(\Omega))\ \Big|\ v_t \in L^2(0,T;L^2(\Omega))\Big\}\,; \\
{\mathcal V}^3(T) =&\ \Big\{ v\in L^2(0,T;H^3(\Omega))\ \Big|\ v_t
\in
L^2(0,T;H^1(\Omega))\Big\}\,; \\
{\mathcal V}^k(T) =&\ \Big\{ v\in L^2(0,T;H^k(\Omega))\ \Big|\ v_t
\in {\mathcal V}^{k-2}(T) \Big\} \quad\text{for $k\ge 4$ }
\end{align*}
with norms
\begin{align*}
\|v\|^2_{{\mathcal V}^1(T)} =&\ \|v\|^2_{L^2(0,T;H^1(\Omega))} + \|v_t\|^2_{L^2(0,T;H^1(\Omega)')} \,;\\
\|v\|^2_{{\mathcal V}^2(T)} =&\ \|v\|^2_{L^2(0,T;H^2(\Omega))} + \|v_t\|^2_{L^2(0,T;L^2(\Omega))} \,;\\
\|v\|^2_{{\mathcal V}^3(T)} =&\ \|v\|^2_{L^2(0,T;H^3(\Omega))} + \|v_t\|^2_{L^2(0,T;H^1(\Omega))} \,;\\
\|v\|^2_{{\mathcal V}^k(T)} =&\ \|v\|^2_{L^2(0,T;H^k(\Omega))} +
\|v_t\|^2_{{\mathcal V}^{k-2}(T)} \quad\text{for $k\ge 4$ }\,.
\end{align*}
We then introduce the space (of ``divergence free'' vector fields)
\begin{align*}
{\mathcal V}_v=\Big\{w\in H^1(\Omega)\ \Big|\ a_i^j(t)w^i_{,j} = 0\
\forall\ t\in[0,T]\Big\}
\end{align*}
and
\begin{align*}
{\mathcal V}_v(T)=\Big\{w\in L^2(0,T;H^1(\Omega))\ \Big|\
a_i^j(t)w^i_{,j} = 0\ \forall\ t\in[0,T]\Big\}\,,
\end{align*}
where the matrix $a$ is defined by (\ref{adefn}). Let $n(\eta) =
(-\eta_2',\eta_1')/|\eta'|$ denote the outward unit normal to $\Gamma(t)$
at the point $\eta(x,t)$.
We define the space $E_\eta^s$ as
\begin{align*}
E_\eta^s = \Big\{ \zeta\in H^2(\Gamma) \ \Big|\
\zeta ''\cdot n \in H^{s-1}(\Gamma)\,,
\zeta'\cdot \eta' \in H^{s-1}(\Gamma) \Big\} \,,
\end{align*}
with norm
\begin{align*}
\|\zeta\|^2_{E_\eta^s} =
\Big[\|\zeta'' \cdot n\|^2_{H^{s-1}(\Gamma)} +
\|\zeta'\cdot \eta' \|^2_{H^{s-1}(\Gamma)}\Big] \,.
\end{align*}

When then set
\begin{align*}
E_\eta^s(T) &= \Big\{ \zeta(t) \in E_\eta^s \ \ t\ a.e.\ \Big|\
\int_0^T \|\zeta(s)\|^2_{E_\eta^s} ds < \infty \Big\}
\end{align*}
with norm
\begin{align*}
\|\zeta\|^2_{E_\eta^s(T)} = \int_0^T \|\zeta(s)\|^2_{E_\eta^s} ds
\,.
\end{align*}

\section{The main theorem}\label{mainthm}
\begin{theorem}\label{maintheorem}
Let $\nu>0$ be given, and
$$F\in L^2(0,T;H^2(\Omega))\cap L^\infty(0,T;H^1(\Omega)), F_t\in L^2(0,T;L^2(\Omega)), F(0)\in H^1(\Omega).$$
Assume that $\Gamma$ is of class $H^{3.5}$ and that the initial data
$u_0\in H^2(\Omega)$ with $\div u_0=0$.  Then
there exists $T>0$ depending on $u_0$ and $F$ such that there exists
a solution $v \in {\mathcal V}^3(T)$ of problem (\ref{solidshelleq})
with $b \in L^2(0,T; H^{2.5}(\Gamma))$ and
$g \in L^2(0,T; H^{2.5}(\Gamma))$.
Moreover, if the initial data has the regularity $u_0 \in H^4(\Omega)$, then
the solution $v \in L^2(0,T;{\mathcal V}^5(\Omega))$ is unique.
\end{theorem}

\section{Preliminary results}
\subsection{Pressure as a Lagrange multiplier}\label{LML}
In the following discussion, we use $\H1H2$ to denote the space
$H^1(\Omega)\cap H^2(\Gamma)$ with norm
$$\|u\|^2_{\H1H2} = \|u\|^2_{H^1(\Omega)} + \|u\|^2_{H^2(\Gamma)}$$
and $\bar{\mathcal V}_\ttv$ ($\bar{\mathcal V}_\ttv(T)$) to denote
the space
\begin{align*}
\Big\{v\in {\mathcal V}_\ttv \ \Big|\ v\in H^2(\Gamma)\Big\} \Big(
\Big\{v\in {\mathcal V}_\ttv(T) \ \Big|\ v\in
L^2(0,T;H^2(\Gamma))\Big\} \Big).
\end{align*}

\begin{lemma}\label{divprob}
For all $p\in L^2(\Omega)$, $t\in [0,T]$, there exists a constant
$C>0$ and $\phi\in \H1H2$ such that $a_i^j(t)\phi_{,j}^i = p$ and
\begin{align}
\|\phi\|_{\H1H2} \le C\|p\|_{L^2(\Omega)}.\label{divestimate}
\end{align}
\end{lemma}
\begin{proof} We solve the following problem on the time-dependent domain $\Omega(t)$:
\begin{align*}
\div (\phi\circ\eta(t)^{-1}) = p\circ\eta(t)^{-1}
\qquad\text{in}\quad\eta(t,\Omega) := \Omega(t).
\end{align*}
The solution to this problem can be written as the sum of the
solutions to the following two problems
\begin{alignat}{2}
\div (\phi\circ\eta(t)^{-1}) =&\ p\circ\eta(t)^{-1} -\bar{p}(t) &&\qquad\text{in}\quad\eta(t,\Omega), \label{divprob1}\\
\div (\phi\circ\eta(t)^{-1}) =&\ \bar{p}(t)
&&\qquad\text{in}\quad\eta(t,\Omega), \label{divprob2}
\end{alignat}
where $\displaystyle{\bar{p}(t)=\frac{1}{|\Omega|}\int_{\Omega}
p(t,x) dx}$. The existence of the solution to problem
(\ref{divprob1}) with zero boundary condition is standard (see, for
example, \cite{GPG1} Chapter 3), and the solution to problem
(\ref{divprob2}) can be chosen as a linear function (linear in $x$)
, for example, $\bar{p}(t) x_1$. The estimate (\ref{divestimate})
follows from the estimates of the solutions to (\ref{divprob1}).
\end{proof}

Define the linear functional on $\H1H2$ by
$(p,a_i^j(t)\varphi_{,j}^i)_{L^2(\Omega)}$ where $\varphi \in
\H1H2$. By the Riesz representation theorem, there is a bounded
linear operator $Q(t): L^2(\Omega)\to \H1H2$ such that for all
$\varphi\in\H1H2$,
\begin{align*}
(p,a_i^j(t)\varphi_{,j}^i)_{L^2(\Omega)} = (Q(t)p,\varphi)_{\H1H2}
:= (Q(t)p,\varphi)_{H^1(\Omega)} + (Q(t)p,\varphi)_{H^2(\Gamma)}.
\end{align*}
Letting $\varphi=Q(t)p$ shows that
$$\|Q(t)p\|_{\H1H2} \le C\|p\|_{L^2(\Omega)}$$
for some constant $C>0$. By Lemma \ref{divprob},
$$\|p\|^2_{L^2(\Omega)}\le \|Q(t)p\|_{\H1H2}\|\varphi\|_{\H1H2} \le C\|Q(t)p\|_{\H1H2}\|p\|_{L^2(\Omega)}$$
which shows that $R(Q(t))$ is closed in $\H1H2$. Since
$\bar{\mathcal V}_v(t)\subset R(Q(t))^\perp$ and $R(Q(t))^\perp
\subset \bar{\mathcal V}_v(t)$, it follows that
\begin{align}
\H1H2(t)=R(Q(t))\oplus_{\H1H2} \bar{\mathcal
V}_v(t).\label{H1decomp}
\end{align}
We can now introduce our Lagrange multiplier
\begin{lemma}\label{lagrangemultiplier}
Let ${\mathcal L}(t)\in \H1H2'$ be such that ${\mathcal L}(t)\varphi
= 0$ for any $\varphi\in \bar{\mathcal V}_v(t)$. Then there exist a
unique $q(t)\in L^2(\Omega)$, which is termed the pressure function,
satisfying
$$\forall\ \varphi\in \H1H2,\quad {\mathcal L}(t)(\varphi)=(q(t),a_i^j(t)\varphi_{,j}^i)_{L^2(\Omega)}.$$
Moreover, there is a $C>0$ (which does not depend on $t\in [0,T]$
and $\epsilon$ and on the choice of $v\in C_T(M)$) such that
$$\|q(t)\|_{L^2(\Omega)} \le C\|{\mathcal L}(t)\|_{\H1H2'}.$$
\end{lemma}
\begin{proof}
By the decomposition (\ref{H1decomp}), for given $\ta$, let
$\varphi=v_1+v_2$, where $v_1\in {\mathcal V}_v(t)$ and $v_2\in
R(Q(t)$. It follows that
$${\mathcal L}(t)(\varphi)={\mathcal L}(t)(v_2) = (\psi(t),v_2)_{\H1H2} = (\psi(t),\varphi)_{\H1H2}$$
for a unique $\psi(t)\in R(Q(t))$.\\
>From the definition of $Q(t)$ we then get the existence of a unique
$q(t)\in L^2(\Omega)$ such that
$$\forall\ \varphi\in \H1H2,\quad{\mathcal L}(t)(\varphi)=(q(t),a_i^j(t) \varphi_{,j}^i)_{L^2(\Omega)}.$$
The estimate stated in the lemma is then a simple consequence of
(\ref{divestimate}).
\end{proof}

\subsection{A polynomial-type inequality}
\begin{lemma}\label{polytypeineq}
Suppose that $x(t)$ is continuous in $[0,T]$, and there are $C_1$,
$C_2$ and $\delta\in (0,1)$ so that
\begin{align*}
x(t) \le C_1 + C_2 t^\delta {\mathcal P}(x(t)) \qquad\forall\
t\in[0,T]\,,
\end{align*}
where ${\mathcal P}$ is a polynomial. Then there $T_1$
(depending only on $C_1$ and $C_2$) such that
\begin{align*}
x(t)\le 2C_1 \qquad\forall\ t\in [0,T_1]\,.
\end{align*}
\end{lemma}
\begin{proof}
We can assume that ${\mathcal P}(x)$ can be factored as $x{\mathcal
Q}(x)$ since the constant part can be collected into $C_1$.
Therefore, we have
\begin{align*}
x(t) \le C_1 + C_2 t^\delta x(t){\mathcal Q}(x(t)), \qquad t\in
[0,T]
\end{align*}
and hence
\begin{align*}
x(t)\Big[1 - C_2 t^\delta {\mathcal Q}(x(t))\Big] \le C_1, \qquad
t\in [0,T].
\end{align*}
Let $T_1>0$ so that $C_2 T_1^\delta {\mathcal Q}(2C_1)\le 1/2$, then
\begin{align*}
\frac{1}{2}x(t) \le C_1, \qquad t\in [0,T_1]\,.
\end{align*}
\end{proof}

\section{Regularized and linearized problem}

Given $\tv \in {\mathcal V}^3(T)$ with the associated $\tg,\tb$ in
$L^2(0,T; H^{2.5}(\Gamma))$,
set $\tF = f \circ\teta$ and
\begin{align*}
\tilde\L_b(\eta) = 2 \Big[|\eta_0'|^{-3} (\eta''\cdot\tn -
b_0)\tn\Big]'' + \Big[|\eta_0'|^{-3} \tg^{-1}\tg' (\teta''\cdot\tn -
b_0) \tn\Big]'
\end{align*}
and
\begin{align*} \tilde\L_m(\eta) = - 4 \Big[(\eta'\cdot\teta' -
|\eta_0'|^2)\teta'\Big]'\,.
\end{align*}
with $\tn(\teta) = (-\teta_2',\teta_1')/|\teta'|$.

The solution $v$ of (\ref{solidshelleq}) is found via a limit as
$\kappa \rightarrow 0$ of the fixed-point of the
map $\tv \mapsto \nv$, where $\nv$ is the solution of the linearized
and $\kappa$-regularized problem:

\begin{subequations}\label{regularizedeq1}
\begin{alignat}{2}
\nv=&\ \neta_t && \text{in }(0,T)\times\Omega\,,\\
\nv^i_t - \nu (\ta_\ell^j \ta_\ell^k \nv_{,k}^i)_{,j} =& - (\ta_i^j \nq)_{,j} + \tF^i && \text{in }(0,T)\times\Omega\,,\\
\ta_i^j \nv_{,j}^i =&\ 0 && \text{in }(0,T)\times\Omega \,,\\
\Big[\nu (\ta_i^k \nv_{,k}^j + \ta_j^k \nv^i_{,k}) -
q\delta_{ij}\Big]\ta_j^\ell N_\ell =&\ \L_m(\teta) +
\tilde\L_b(\eta) + \kappa \eta^{(4)} \quad&& \text{on }(0,T)\times\Gamma \,,\\
\nv(0,x)=&\ u_0(x) && \text{on } \{t=0\}\times\Omega\,,\\
\neta =&\ \id && \text{on } \{t=0\}\times\Omega\,,
\end{alignat}
\end{subequations}
where we set $\varepsilon = 1$ and ignore the factor $1/3$ in front
of $\L_b$. Note that here we treat the membrane traction as an extra
forcing on the boundary. Also note that the time $T$ a priori depends on
$\kappa$.

Following the same analysis as in \cite{ChCoSh2006}, we can show that for this
regularized problem
(for a given and fixed $\ttv$), there exists a unique solution $(\neta,\nv)$
to (\ref{regularizedeq1}) with
$\nv \in {\mathcal V}^3(T)$ and $\neta \in L^2(0,T; H^{5.5}(\Gamma))$.

This follows by first
approximating by a penalized problem, and then performing a
regularity analysis (energy estimates). By the Tychonoff fixed-point
theorem, there exists a fixed point $\nv$ in ${\mathcal
V}^3(T_\kappa)$ with $\int_0^t \nv ds \in L^2(0,T_\kappa; H^{5.5}(\Gamma))$ and
$(\na)_i^j v_{\kappa,j}^i = 0$, and this $\nv$ and the associated
$\neta$ satisfy
\begin{subequations}\label{regularizedeq}
\begin{alignat}{2}
\nv=&\ \neta_t && \text{in }(0,T_\kappa)\times\Omega\,,\\
\nv^i_t - \nu [(\na)_\ell^j (\na)_\ell^k \nv_{,k}^i]_{,j} =& - (\na)_i^j \nq_{,j} + F^i && \text{in }(0,T_\kappa)\times\Omega\,,\\
(\na)_i^j v_{\kappa,j}^i =&\ 0&& \text{in }(0,T_\kappa)\times\Omega\,,\\
\Big[\nu (\na_i^k v_{,k}^j + \na_j^k \nv^i_{,k}) -
\nq\delta_{ij}\Big] (\na)_j^\ell N_\ell =&\ \L_m(\neta) + \L_b(\neta) \ && \text{on }(0,T_\kappa)\times\Gamma \,,\label{bc}\\
& + \kappa \neta^{(4)} && \nonumber\\
\nv(0,x)=&\ u_0(x) && \text{on }\{t=0\}\times\Omega\,,\\
\neta =&\ \id && \text{on }\{t=0\}\times\Omega\,.
\end{alignat}
\end{subequations}

%\begin{remark} The reason of needing the regularization in $\eta$ is that
%the normal vector $\tn$ constructed from the previous step $\teta$
%(by the formula $\tn = \tg^{1/2} \teta_{,1}\times \teta_{,2}$) is
%always less regular than $\teta$. Therefore, we add artificial
%viscosity into the boundary condition so that the solution to this
%regularized problem produces a much more regular normal and hence
%the solution $v$ is in the same space as the initial input $\tv$.
%However, the price we pay is that the fixed point $\nv$ only exists
%in ${\mathcal V}^3(T_\kappa)$ with a $\kappa$-dependent time period
%$T_\kappa$. In general, $T_\kappa \to 0$ as $\kappa\to 0$ and the
%limit of $\nv$ as $\kappa\to 0$ means nothing to the problem. See
%Remark \ref{needH35normal} for more details.
%\end{remark}

%\begin{remark} Unlike the case in \cite{ChCoSh2006}, we cannot obtain a
%$\kappa$-independent estimate for the regularized problem.
%Therefore, we keep this artificial viscosity at the level of
%applying the Tychonoff fixed-point theorem to obtain a fixed point
%%to (\ref{regularizedeq1}). We then only need to find a
%$\kappa$-independent estimate so that $T_k$ can be chosen
%independent of $\kappa$. The limit of $u_\kappa$ then will be a
%solution to (\ref{solidshelleq}).
%\end{remark}

\section{A priori estimates for $\nv$, $\nq$ and $\neta$}\label{knownfacts}
%\subsection{Regularity of $\nv$, $\nq$ and $\neta$}

\subsection{$L^2(\Omega)$-estimate for $\nq$} By (\ref{regularizedeq}), a solution $\nv$, $\nq$ and
$\neta$ also satisfy
\begin{align}
& \langle \nv_t,\phi\rangle + \frac{\nu}{2} \langle (\na)_i^k
(v_\kappa^i)_{,k} + (\na)_j^k (v_\kappa^j)_{,k}, (\na)_i^k
\phi_{,k}^i + (\na)_j^k \phi_{,k}^j\rangle + \langle \nq, (\na)_i^j\phi_{,j}^i\rangle \nonumber\\
+& \int_\Gamma |\eta_0'|^{-3}\Big[ 4 (|\eta_\kappa'|^2 -
|\eta_0'|^2) (\eta_\kappa'\cdot \phi') + 2 (\eta_\kappa''\cdot \nn - b_0) (\phi''\cdot \nn)\Big] dS \label{weakform}\\
-& \int_\Gamma |\eta_0'|^{-3} g_\kappa^{-1} (\eta_\kappa''\cdot \nn
- b_0) g_\kappa' (\phi'\cdot \nn) dS + \kappa \int_\Gamma
\eta_\kappa'' \cdot \phi'' dS= \langle F,\phi\rangle \nonumber
\end{align}
for all $\phi\in\H1H2$. Therefore, by the Lagrange multiplier lemma,
\begin{align}
& \|\nq\|^2_{L^2(\Omega)} \le C\Big[\|\nv_t\|^2_{L^2(\Omega)} +
\|D_\neta \nv\|^2_{L^2(\Omega)} + \|(\ng -
g_0)\eta_\kappa'\|^2_{H^{-1}(\Gamma)} \nonumber \\
&\qquad + \|(\nb - b_0) \nn\|^2_{L^2(\Gamma)} + \|g_\kappa^{-1} (\nb
- b_0) g_\kappa'\nn\|^2_{H^{-1}(\Gamma)} + \kappa
\|\neta\|^2_{H^2(\Gamma)}
\Big] \nonumber \\
\le&\ C \mP \Big[\|\nv_t\|^2_{L^2(\Omega)} +
\|\nabla\nv\|^2_{L^2(\Omega)} + \|\ng - g_0\|^2_{L^2(\Gamma)} \label{qL2estimate}\\
&\qquad\qquad\qquad\quad\ \ + \|\nb - b_0\|^2_{L^2(\Gamma)} + \kappa
\|\eta\|^2_{H^2(\Gamma)}\Big] \nonumber
\end{align}
for some constant $C$ independent of $\kappa$.

\subsection{Interior regularity} Converting the fluid equation
(\ref{regularizedeq}b) into Eulerian variables by composing with
$\eta_\kappa^{-1}$, we obtain a Stokes problem in the domain
$\neta(\Omega)$:
\begin{subequations}\label{stokes}
\begin{align}
-\nu\Delta u_\kappa + \nabla p_\kappa &= F\circ\eta_\kappa^{-1} -
\nv_t\circ\eta_\kappa^{-1} %+ \nu\na_{\ell,j}^j\circ\eta_\kappa^{-1}
%u_{\kappa,\ell} - p_\kappa \na^j_{i,j}\circ\eta_\kappa^{-1}
\,, \\
\div u_\kappa &= 0\,,
\end{align}
\end{subequations}
where $u_\kappa = \nv\circ \eta_\kappa^{-1}$ and $p_\kappa =
\nq\circ \eta_\kappa^{-1}$. By the regularity results for the Stokes
problem,
\begin{align*}
& \|u_\kappa\|^2_{H^2(\neta(\Omega))} +
\|p_\kappa\|^2_{H^1(\neta(\Omega))} \\ \le&\ C
\Big[\|F\circ\eta_\kappa^{-1}\|^2_{L^2(\neta(\Omega))} +
\|\nv_t\circ\eta_\kappa^{-1}\|^2_{L^2(\neta(\Omega))}%\Big] \\
%&\qquad\qquad + C \mP \Big[\|\nabla
%u_\kappa\|^2_{L^2(\neta(\Omega))} +
%\|p_\kappa\|^2_{L^2(\neta(\Omega))}
+ \|u_\kappa\|^2_{H^{1.5}(\Gamma)} \Big]
\end{align*}
or
%\begin{align*}
%\|\nv\|^2_{H^2(\Omega)} + \|\nq\|^2_{H^1(\Omega)} \le&\
%C\mP\Big[\|F\|^2_{L^2(\Omega)} + \|\nv_t\|^2_{L^2(\Omega)} +
%\|\nabla \nv\|^2_{L^2(\Omega)} \\
%&\qquad\qquad\qquad\quad\ \ + \|\nv\|^2_{H^{1.5}(\Gamma)} +
%\|\nq\|^2_{L^2(\Omega)} \Big] \,.
%\end{align*}
%By the Lagrange multiplier lemma,
\begin{align}
\|\nv\|^2_{H^2(\Omega)} + \|\nq\|^2_{H^1(\Omega)} \le C \Big[
\|F\|^2_{L^2(\Omega)} + \|\nv_t\|^2_{L^2(\Omega)} \Big] + C \mP
\|\nv\|^2_{H^{1.5}(\Gamma)} \,. \label{regularity1temp}
\end{align}
Similarly,
%\begin{align*}
%\|\nv\|^2_{H^3(\Omega)} + \|\nq\|^2_{H^2(\Omega)} \le&\ C\mP\Big[
%\|F\|^2_{H^1(\Omega)} + \|\nv_t\|^2_{H^1(\Omega)} + \|\nabla
%\nv\|^2_{H^1(\Omega)} \\
%&\qquad\qquad\qquad\quad\ \ + \|\nv\|^2_{H^{2.5}(\Gamma)} +
%\|\nq\|^2_{H^1(\Omega)} \Big]
%\end{align*}
%and therefore by (\ref{regularity1temp}),
\begin{align}
& \|\nv\|^2_{H^3(\Omega)} + \|\nq\|^2_{H^2(\Omega)} \le
C\mP\Big[\|F\|^2_{H^1(\Omega)} + \|\nv_t\|^2_{H^1(\Omega)} +
\|\nv\|^2_{H^{2.5}(\Gamma)} \Big]\,. \label{regularity2temp}
\end{align}

\subsection{$H^1(\Omega)$-estimate for $\nv_t$}
\begin{align}
\|\nabla \nv\|^2_{L^2(\Omega_0)} =&\ \|\nabla u_0 + \int_0^t \nabla
\nv_t ds\|^2_{L^2(\Omega_0)}
\le \Big[\|\nabla u_0\|_{L^2(\Omega_0)} + \int_0^t \|\nabla \nv_t\|_{L^2(\Omega_0)} ds \Big]^2 \nonumber\\
\le&\ 2\Big[\|u_0\|^2_{H^1(\Omega_0)} + t
\|\nv_t\|^2_{L^2(0,T;H^1(\Omega))}\Big] \label{H1estimate}\,.
\end{align}

\section{Elliptic estimates on the boundary}
\subsection {Estimates without artificial viscosity}
Since $\nv\in {\mathcal V}^3(T)$, the associated $\neta$ satisfies
the boundary condition (\ref{bc}) in the pointwise sense. We start
with the estimates without considering the artificial viscosity to
illustrate the basic idea; then in the next section we consider the
full boundary condition (\ref{bc}) and obtain the desired estimates.
By (\ref{normalderivatives}), we find that
\begin{align*}
\L_m(\eta) = 4 \Big[\sqrt{g} (g' - g_0') +
\frac{g-g_0}{2\sqrt{g}}g'\Big] \tau + 4(g - g_0) bn
\end{align*}
and
\begin{align*}
\L_b(\eta) =&\ \Big[2 \Big(|\eta_0'|^{-3}(b-b_0)\Big)'' - 2
|\eta_0'|^{-3} g^{-1} b^2 (b-b_0) +
\Big(|\eta_0'|^{-3} g^{-1} g'(b-b_0) \Big)'\Big] n \\
& - \Big[4\Big(|\eta_0'|^{-3}(b-b_0)\Big)' + 2 |\eta_0'|^{-3}
g^{-1/2} (b-b_0) b'\Big] g^{-1/2} b \tau \,.
\end{align*}
Given $h\in L^2(0,T;H^{1.5}(\Gamma))\cap
L^\infty(0,T;H^{0.5}(\Gamma))$, a solution to
\begin{align}
\L_m(\eta) + \L_b(\eta) = h \label{tempeq}
\end{align}
satisfies the ``normal equation''
\begin{align}
& 2 \Big[|\eta_0'|^{-3}(b-b_0)\Big]'' \label{normaleq}\\
=&\ h\cdot n - 4(g - g_0) b + 2 |\eta_0'|^{-3} g^{-1} b^2 (b-b_0) -
\Big[|\eta_0'|^{-3} g^{-1} g'(b-b_0) \Big]' \,. \nonumber
\end{align}
We also have the ``tangential equation''
\begin{align}
& 4 \sqrt{g} (g - g_0)' \label{tangentialeq}\\
=&\ h\cdot \tau + g^{-1/2} \Big[4\Big(|\eta_0'|^{-3}(b-b_0)\Big)'b +
2 |\eta_0'|^{-3} g^{-1/2} (b-b_0) b b' - 2(g-g_0) g' \Big] \,.
\nonumber
\end{align}
Therefore, by elliptic estimates, a solution to (\ref{tempeq})
satisfies
\begin{subequations}\label{bgdiffestimate}
\begin{align}
\|b-b_0\|^2_{H^{2.5}(\Gamma)} \le&\ C\Big[ \|h\cdot
n\|^2_{H^{0.5}(\Gamma)} + \|(g-g_0) b\|^2_{H^{0.5}(\Gamma)} \\
&\quad + \|g^{-1} b^2 (b-b_0)\|^2_{H^{0.5}(\Gamma)} + \|[g^{-1} g'
(b-b_0)]'\|^2_{H^{0.5}(\Gamma)} \Big]\,, \nonumber \\
\|g-g_0\|^2_{H^{2.5}(\Gamma)} \le&\ C {\mathcal
P}(\|\eta\|^2_{H^{2.5}(\Gamma)}) \Big[ \|h\cdot
\tau\|^2_{H^{1.5}(\Gamma)} + \|(b-b_0)'g^{-1/2} b
\|^2_{H^{1.5}(\Gamma)} \\
&\qquad\quad + \|g^{-1} (b-b_0) b'b\|^2_{H^{1.5}(\Gamma)} + \|(g -
g_0)g'\|^2_{H^{0.5}(\Gamma)} \Big] \,,\nonumber
\end{align}
\end{subequations}
where $C$ only depends on $\Gamma$. Since
\begin{subequations}\label{interpolation}
\begin{align}
\|fg\|_{H^{0.5}(\Gamma)} &\le C \Big[\|f\|_{H^{0.5}(\Gamma)}
\|g\|_{L^\infty(\Gamma)} + \|f\|_{H^{0.5}(\Gamma)}
\|g\|_{L^\infty(\Gamma)}\Big]\,, \\
\|fg\|_{H^{1.5}(\Gamma)} &\le C \Big[\|f\|_{H^{1.25}(\Gamma)}
\|g\|_{H^{1.5}(\Gamma)} + \|f\|_{H^{1.5}(\Gamma)}
\|g\|_{H^{1.25}(\Gamma)}\Big]\,,
\end{align}
\end{subequations}
by the Leibnitz rule, we have
\begin{align}
& \|[g^{-1} g' (b-b_0)]'\|^2_{H^{0.5}(\Gamma)} \nonumber\\
\le&\ C {\mathcal P}(\|\eta\|^2_{H^{2.5}(\Gamma)})
\Big[\|g\|^2_{H^{2.5}(\Gamma)} \|b-b_0\|^2_{L^\infty(\Gamma)} + \|b
- b_0\|^2_{W^{1,\infty}(\Gamma)} \Big] \nonumber\\
\le&\ C_\epsilon {\mathcal P}(\|\eta\|^2_{H^{2.5}(\Gamma)})
\Big[\|g\|^2_{H^{2.5}(\Gamma)} \|b-b_0\|^2_{L^\infty(\Gamma)} +
\|b-b_0\|^2_{H^{1.5}(\Gamma)}\Big] \label{btempineq}\\
& + \epsilon \|b-b_0\|^2_{H^{2.5}(\Gamma)} \,, \nonumber
\end{align}
and
\begin{align*}
\|g^{-1} (b-b_0) b'b\|^2_{H^{1.5}(\Gamma)} \le&\ C {\mathcal
P}(\|\eta\|^2_{H^{2.5}(\Gamma)}) \|b-b_0\|^2_{H^{1.25}(\Gamma)}
\|b\|^2_{H^{2.5}(\Gamma)} \|b\|^2_{H^{1.5}(\Gamma)}\,.
\end{align*}
Let $X(T) = \|v\|^2_{{\mathcal V}^3(T)} +
\|b\|^2_{L^2(0,T;H^{2.5}(\Gamma))} +
\|g\|^2_{L^2(0,T;H^{2.5}(\Gamma))}$, then
\begin{subequations}\label{bgsmall}
\begin{align}
\|b(t) - b_0\|^2_{H^{1.25}(\Gamma)} +
\|g(t)-g_0\|^2_{H^{1.25}(\Gamma)} &\le C \sqrt{t} {\mathcal P}(X(T)) \,, \\
%\le C \sqrt{t} {\mathcal P}(\|\eta\|^2_{H^{2.5}(\Gamma)}) \Big[
%\|b-b_0\|^2_{L^2(0,T;H^{2.5}(\Gamma))} &+
%\|g-g_0\|^2_{L^2(0,T;H^{2.5}(\Gamma))}\Big] \\
%\qquad\quad \times \Big[\|v_0\|^2_{H^2(\Omega)} +
%\|v\|^2_{{\mathcal V}(T)} \Big]\\
\|b(t) - b_0\|^2_{H^{1.5}(\Gamma)} +
\|g(t)-g_0\|^2_{H^{1.5}(\Gamma)}&\le C {\mathcal P}(X(T))\,,
\end{align}
\end{subequations}
for all $0\le t\le T$. Therefore, by choosing $\epsilon>0$ small
enough in (\ref{btempineq}),
\begin{align}
\|b\|^2_{H^{2.5}(\Gamma)} %\le &\ C {\mathcal
%P}(\|\eta\|^2_{H^{2.5}(\Gamma)}) \Big[\|f\|^2_{H^{0.5}(\Gamma)}
%\|n\|^2_{L^\infty(\Gamma)} + %\|b\|_{H^{0.5}(\Gamma)}
%\|b\|^2_{L^\infty(\Gamma)}
%\|b-b_0\|^2_{L^\infty(\Gamma)} \\
%&\qquad\qquad + \|g\|^2_{H^{2.5}(\Gamma)}
%\|b-b_0\|^2_{L^\infty(\Gamma)} + \|g-g_0\|^2_{L^\infty(\Gamma)}
%\|b\|^2_{H^{0.5}(\Gamma)} \\
%&\qquad\qquad + \|b-b_0\|^2_{H^{1.5}(\Gamma)} \Big] \\
%\le&\ C {\mathcal P}(\|\eta\|^2_{H^{2.5}(\Gamma)})
%\Big[\|f\|^2_{H^{1.5}(\Gamma)} + \|f\|^2_{H^{0.5}(\Gamma)}
%\|b\|^2_{H^{1.5}(\Gamma)} \\
%& + \Big( \|b-b_0\|^2_{H^{1.25}(\Gamma)} +
%\|g-g_0\|^2_{H^{1.25}(\Gamma)} \Big)\Big( \|b\|^2_{H^{2.5}(\Gamma)}
%+ \|g\|^2_{H^{2.5}(\Gamma)}\Big) \Big] \\
\le&\ C {\mathcal P}(\|\eta\|^2_{H^{2.5}(\Gamma))})
\|h\|^2_{L^\infty(\Gamma)}+ C {\mathcal P}(X(T)) \Big[1 +
\|h\|^2_{H^{0.5}(\Gamma)}\Big] \label{bestimate}\\
& + C \sqrt{t} {\mathcal P}(X(T)) \Big[\|b\|^2_{H^{2.5}(\Gamma)} +
\|g\|^2_{H^{2.5}(\Gamma)} \Big] \nonumber
\end{align}
and hence by (\ref{bgdiffestimate}a) (and also (\ref{bgsmall}),
(\ref{bestimate})),
\begin{align}
\|g\|^2_{H^{2.5}(\Gamma)} \le&\ C {\mathcal
P}(\|\eta\|^2_{H^{2.5}(\Gamma)}) \|h\|^2_{H^{1.5}(\Gamma)} + C%_\epsilon
P(X(T)) \Big[1 + \|h\|^2_{H^{0.5}(\Gamma)}\Big]
%+ \epsilon \|b-b_0\|^2_{H^{2.5}(\Gamma)}
\label{gestimate}\\
& + C \sqrt{t} {\mathcal P}(X(T)) \Big[\|b\|^2_{H^{2.5}(\Gamma)} +
\|g\|^2_{H^{2.5}(\Gamma)}\Big] \,. \nonumber
\end{align}
With $h = \Big[\mu (a_i^k v_{,k}^j + a_j^k
v^i_{,k})-q\delta_{ij}\Big]a_j^\ell N_\ell$ in mind, we find that
\begin{align}
\int_0^t \|b(s)\|^2_{H^{2.5}(\Gamma)} ds \le&\ C (t+ t^{1/2} +
t^{1/4}) {\mathcal P}(X(T)) \label{bL2H25estimate}
\end{align}
and
\begin{align}
\int_0^t \|g(s)\|^2_{H^{2.5}(\Gamma)} ds \le&\ C (1 + t^{1/2} + t)
{\mathcal P}(X(T)) \,.\label{gL2H25estimate}
\end{align}

%\begin{remark}
%It appears that (\ref{tempeq}) is elliptic only when $T$ is small
%enough when $n=2$. It is even not elliptic when $n=3$. See Section
%\ref{3dfail} for more details.
%\end{remark}

\subsection{Estimates with artificial viscosity}
Now we study the full boundary condition
\begin{align}
\L_m(\neta) + \L_b(\neta) + \kappa \eta_\kappa''''= h \,.
\label{tempeq1}
\end{align}
By the Leibnitz rule and (\ref{normalderivatives}),
\begin{subequations}\label{eta4prime}
\begin{align}
\eta_\kappa''''\cdot \nn =&\ b_\kappa'' - g_\kappa^{-1} b_\kappa^3 -
\frac{3}{4} g_\kappa^{-2} (g_\kappa')^2 \nb + \frac{1}{2}
g_\kappa^{-1} g_\kappa' b_\kappa' + g_\kappa^{-1} g_\kappa'' \nb\,, \\
\eta_\kappa''''\cdot \tau =&\ \frac{1}{2} g_\kappa^{-1/2} \Big[
g_\kappa''' - 3 \nb b_\kappa' - \frac{3}{2} g_\kappa^{-1}
g_\kappa'g_\kappa'' + \frac{3}{4} g_\kappa^{-2} (g_\kappa')^3\Big]
\,,
\end{align}
\end{subequations}
where $\ng = \eta_\kappa'\cdot \eta_\kappa'$, $\nn = g_\kappa^{-1/2}
\neta_{,1}\times \neta_{,2}$ and $\nb = \eta_\kappa''\cdot \nn$.
Define
\begin{align*}
\nX(T) =& \sup_{0\le t\le T} \Big[\|\nv_t\|^2_{L^2(\Omega)} +
\|\nv\|^2_{H^2(\Omega)} + \|\ng\|^2_{H^2(\Gamma)} +
\|\nb\|^2_{H^2(\Gamma)}
+ \|v_\kappa'\cdot \eta_\kappa'\|^2_{L^2(\Gamma)} \nonumber\\
&\quad + \|v_\kappa''\cdot \nn\|^2_{L^2(\Gamma)} \Big] +
\|\nv\|^2_{{\mathcal V}^3(T)} + \|\nb\|^2_{L^2(0,T;H^{2.5}(\Gamma))}
+ \|\ng\|^2_{L^2(0,T;H^{2.5}(\Gamma))}\,.
\end{align*}
By the same technique, we find using (\ref{normaleq}) that
\begin{align*}
& \|\nb\|^2_{H^{2.5}(\Gamma)} \le C {\mathcal
P}(\|\neta\|^2_{H^{2.5}(\Gamma)}) \|h\|^2_{L^\infty(\Gamma)} + C
{\mathcal
P}(\nX(T))\Big[1+\|h\|^2_{H^{0.5}(\Gamma)}\Big] \\
&\qquad + C\sqrt{t} P(\nX(T)) \Big[\|\nb\|^2_{H^{2.5}(\Gamma)} +
\|\ng\|^2_{H^{2.5}(\Gamma)} \Big] + C_\epsilon \kappa {\mathcal
P}(\nX(T)) \\
&\qquad + \frac{\epsilon \kappa}{{\mathcal P}(\nX(T))}
\|\ng\|^2_{H^3(\Gamma)}
\end{align*}
and from the tangential equation (\ref{tangentialeq}), we find that
\begin{align*}
& \|\ng\|^2_{H^{2.5}(\Gamma)} + \kappa \|\ng\|^2_{H^{4.5}(\Gamma)} \nonumber\\
\le&\ C {\mathcal P}(\|\neta\|^2_{H^{2.5}(\Gamma)})
\|h\|^2_{H^{1.5}(\Gamma)} + C {\mathcal
P}(\nX(T))\Big[1+\|h\|^2_{H^{0.5}(\Gamma)}\Big] \\
& + C\sqrt{t} P(\nX(T)) \Big[\|\nb\|^2_{H^{2.5}(\Gamma)} +
\|\ng\|^2_{H^{2.5}(\Gamma)}\Big] + C_\epsilon \kappa {\mathcal
P}(\nX(T)) + C \epsilon \kappa \|\ng\|^2_{H^{4.5}(\Gamma)} \,,
\end{align*}
where we use the fact that
\begin{align*}
&\|\eta_\kappa''''\cdot \nn\|^2_{H^{0.5}(\Gamma)} \\
\le&\ C \kappa \Big[\|\nb\|^4_{L^\infty(\Gamma)}
\|\nb\|^2_{H^{0.5}(\Gamma)} + \|\nb\|^2_{L^\infty(\Gamma)}
\|g_\kappa\|^2_{H^{2.5}(\Gamma)} +
\|\nb\|^2_{H^{1.5}(\Gamma)} \|g_\kappa'\|^2_{L^\infty(\Gamma)} \Big] \\
\le&\ C_\epsilon \kappa {\mathcal P}(\nX(T)) + \epsilon \kappa
\frac{1}{{\mathcal P}(\nX(T))} \|\eta_\kappa'\|^2_{H^3(\Gamma)} \,,
\end{align*}
and
\begin{align*}
&\|\eta_\kappa''''\cdot \tau_\kappa\|^2_{H^{0.5}(\Gamma)} \\
\le&\ C \kappa \Big[\|(\nb-b_0)b_\kappa' + b_0
b_\kappa'\|^2_{H^{1.5}(\Gamma)} + \|\ng\|^2_{H^{2.5}(\Gamma)}
\|\ng\|^2_{H^{3.25}(\Gamma)} +
\|\ng\|^4_{H^{2.25}(\Gamma)} \|\ng\|^2_{H^{2.5}(\Gamma)} \Big] \\
\le&\ C\sqrt{t} P(\nX(T)) \|\nb\|^2_{H^{2.5}(\Gamma)} + C\kappa
\|\nb\|^2_{H^{2.5}(\Gamma)} + C_\epsilon \kappa {\mathcal P}(\nX(T))
+ \epsilon \kappa \|\ng\|^2_{H^{4.5}(\Gamma)} \,.
\end{align*}
Note that here $C$and $C_\epsilon$ are independent of $\kappa$.
Choosing $\epsilon>0$ small enough, we find that
\begin{align}
\Big(\int_0^t \|\nb(s)\|^2_{H^{2.5}(\Gamma)} ds \le\Big)& {\mathcal
P}(\nX(T)) \int_0^t \|\nb(s)\|^2_{H^{2.5}(\Gamma)} ds
\label{bkestimate}\\
\le&\ C (t + t^{1/2} + t^{1/4}) P(\nX(T)) + \frac{\kappa}{100}
\int_0^t \|\eta_\kappa'(s)\|^2_{H^3(\Gamma)} ds \nonumber
\end{align}
and
\begin{align}
\int_0^t \Big[\|\ng(s)\|^2_{H^{2.5}(\Gamma)} + \kappa
\|\ng(s)\|^2_{H^{4.5}(\Gamma)}\Big]ds \le C (1 + t^{1/2} + t)
{\mathcal P}(\nX(T))\,. \label{gkestimate}
\end{align}

\subsection{The estimate of $\nn$} By (\ref{normalderivatives}),
\begin{subequations}\label{normalestimates}
\begin{align}
\|\nn\|^2_{H^{2.5}(\Gamma)} \le&\ %C {\mathcal
%P}(\|\eta\|^2_{H^{2.5}(\Gamma)}) \|b\|^2_{H^{1.5}(\Gamma)} \le
C{\mathcal P}(\nX(T))\,, \\
\int_0^t \|\nn\|^2_{H^{3.5}(\Gamma)} ds \le&\ C (1 + t^{1/2} + t)
{\mathcal P}(\nX(T)) ds\,.
%C {\mathcal P}(\|\eta\|^2_{H^{2.5}(\Gamma)})
%\Big[\|b\|^2_{H^{1.25}(\Gamma)} (\|g\|^2_{H^{2.5}(\Gamma)} +
%\|\eta\|^2_{H^{3.5}(\Gamma)}) +
%\|b\|^2_{H^{2.5}(\Gamma)}\Big] \nonumber\\
%\le&\ C {\mathcal P}(\|\eta\|^2_{H^{2.5}(\Gamma)})
%\|f\|^2_{H^{1.5}(\Gamma)} + C {\mathcal P}(X(T))
%\Big[1+\|f\|^2_{H^{0.5}(\Gamma)}\Big] \\
%&\ + C \sqrt{t} {\mathcal P}(X(T)) \Big[\|b\|^2_{H^{2.5}(\Gamma)} +
%\|g\|^2_{H^{2.5}(\Gamma)}\Big] \,. \nonumber
\end{align}
\end{subequations}

\subsection{Small time results}
In this section, we rewrite some inequalities in Section
\ref{knownfacts} that will be used in the later discussion. First of
all, note that (\ref{H1estimate}) implies that
\begin{align}
\|v\|^2_{H^1(\Omega)} \le C \Big[\|u_0\|^2_{H^1(\Omega)} + t
{\mathcal P}(\nX(T))\Big]\,. \label{H1estimate1}
\end{align}
Since
\begin{align*}
\|\neta\|^2_{H^2(\Gamma)} \le 2 \Big[\|\id\|^2_{H^2(\Gamma)} + t
\int_0^t \|\nv\|^2_{H^{2.5}(\Omega)} ds\Big]\,,
\end{align*}
(\ref{qL2estimate}) can be rewritten as
\begin{align}
\|\nq\|^2_{L^2(\Omega)} \le C \Big[\|\nv_t\|^2_{L^2(\Omega)} + t
{\mathcal P}(\nX(T))\Big]\,. \label{qL2estimate1}
\end{align}
Finally, we also have
\begin{align}
\int_0^t \Big[\|\nv\|^2_{H^2(\Omega)} +
\|\nq\|^2_{H^1(\Omega)}\Big]ds \le&\ C \int_0^t
\|F\|^2_{L^2(\Omega)} ds + C t \nX(T) \label{regularity1}
% \\ & + \epsilon \int_0^t \|\nv\|^2_{H^{2.5}(\Gamma)} ds \nonumber
\end{align}
and
\begin{align}
\int_0^t \Big[\|\nv\|^2_{H^3(\Omega)} +
\|\nq\|^2_{H^2(\Omega)}\Big]ds \le&\ C \int_0^t
\Big[\|F\|^2_{H^1(\Omega)} + \|\nv_t\|^2_{H^1(\Omega)} +
\|\nv\|^2_{H^{2.5}(\Gamma)}\Big] ds \nonumber\\
& + C (t + \sqrt{t}) {\mathcal P}(\nX(T)) \,,\label{regularity2}
\end{align}
where we have used that
\begin{align}
\|\nabla \neta(t) - \id \|_{H^2(\Omega)} + \|\na(t) - \id \|_{H^2(\Omega)} %&\le C \sqrt{t} \nX(T)(1 + \sqrt{t}
%\nX(T))
&\le C (t + \sqrt{t}) {\mathcal P}(\nX(T)) \label{aetaineq}
\end{align}
to remove the $\na$ and $\neta$ dependence from inequalities
(\ref{regularity1temp}) and (\ref{regularity2temp}).

\section{Nonlinear estimates}
In the following discussion, we will always assume that $T\le 1$.
Therefore, all the time dependent functions appearing in the
previous section, such as $t$, $(t + \sqrt{t})$, $(t+t^{1/2} +
t^{1/4})$, etc., can be replaced by $t^\delta$ for some fixed
$\delta \in (0,1)$.
\subsection{Partition of unity}\label{Kestimate} Since $\Omega$ is compact, by partition of unity, we can choose two non-negative smooth
functions $\zeta_0$ and $\zeta_1$ so that
\begin{align*}
\zeta_0 + \zeta_1 = 1 \quad\text{in}\quad \Omega \,;\quad
\supp(\zeta_0) \subset\subset \Omega \,;\quad \supp(\zeta_1)
\subset\subset \Gamma\times(-\epsilon,\epsilon):=\Omega_1.
\end{align*}
We will assume that $\zeta_1=1$ inside the region $\Omega_1'\subset
\Omega_1$ and $\zeta_0=1$ inside the region $\Omega'\subset \Omega$.
Note that then $\zeta_1 = 1$ while $\zeta_0 = 0$ on $\Gamma$.
\subsection{Energy estimates for $v_\kappa''$ near the boundary} Use
$(\zeta_1^2 v_\kappa'')''$ as a test function in (\ref{weakform}),
we find that
\begin{align}
&\ \ \frac{1}{2}\frac{d}{dt} \Big[\|\zeta_1
v_\kappa''\|^2_{L^2(\Omega)} + 2 \int_\Gamma |\eta_0'|^{-3} [|(\ng -
g_0)''|^2 + |(\nb -
b_0)''|^2 ] dS + \kappa \int_\Gamma |\eta_\kappa^{(4)}|^2 dS\Big] \nonumber\\
&+ \frac{\nu}{2} \|\zeta_1 D_\neta v_\kappa''\|^2_{L^2(\Omega)} =
\langle F, (\zeta_1^2 v_\kappa'')'' \rangle + (K_1 + K_2 + \cdots +
K_7 )\,, \label{L2H3eq}
\end{align} where $(D_\neta w)_i^j \equiv (\na)_i^k w_{,k}^j + (\na)_j^k
w_{,k}^i$ is the nonlinear version of the rate of deformation
tensor, and $K_i's$ are defined by
\begin{align*}
K_1 =&\ - 2\nu\int_{\Omega} \zeta_1 \Big[(a_\kappa')_i^k
(v_\kappa^j)_{,k}' + (a_\kappa')_j^k (v_\kappa^i)_{,k}'\Big]
(\na)_i^k
\zeta_1 (v_\kappa^j)_{,k}'' dx \\
& - \nu \int_{\Omega} \zeta_1 \Big[(a_\kappa'')_i^k
(v_\kappa^j)_{,k} + (a_\kappa'')_j^k (v_\kappa^i)_{,k}\Big]
(\na)_i^k \zeta_1 (v_\kappa^j)_{,k}'' dx \\
& - 2 \nu \int_{\Omega} \zeta_1 \Big[(a_\kappa')_i^k
(v_\kappa^j)_{,k} + (a_\kappa')_j^k (v_\kappa^i)_{,k}\Big]
(a_\kappa')_i^k \zeta_1 (v_\kappa^j)_{,k}'' dx \\
& - 2 \nu \int_{\Omega} \zeta_1 \Big[(\na)_i^k (v_\kappa^j)_{,k}' +
(\na)_j^k (v_\kappa^i)_{,k}'\Big] (a_\kappa')_i^k \zeta_1
(v_\kappa^j)_{,k}'' dx\,,
\end{align*}
\begin{align*}
K_2 =&\ 4\int_\Gamma |\eta_0'|^{-3} (\ng - g_0)
\Big[\eta_\kappa^{(5)}\cdot v_\kappa' + 4 \eta_\kappa^{(4)}\cdot
v_\kappa'' + 6 \eta_\kappa'''\cdot v_\kappa''' + 4\eta_\kappa''\cdot
v_\kappa^{(4)}\Big] dS\,,\\
K_3 =&\ - 8 \int_\Gamma (\|\eta_0'|^{-3})'(\ng-g_0)'
(\eta_\kappa'\cdot v_\kappa')'' dS - 4 \int_\Gamma
(|\eta_0'|^{-3})''(\ng-g_0) (\eta_\kappa'\cdot v_\kappa')'' dS\\
&\ - 4\int_\Gamma \Big[|\eta_0'|^{-3} (g_\kappa - g_0)\Big]'' g_0''
dS\,,\\
K_4 =&\ 2 \int_\Gamma |\eta_0'|^{-3} (\nb - b_0) \Big[
(\eta_\kappa''\cdot n_{\kappa t})^{(4)} + n_\kappa^{(4)}\cdot
v_\kappa'' + 4 n_\kappa'''\cdot v_\kappa''' + 6 n_\kappa''\cdot v_\kappa^{(4)} \\
&\qquad\qquad\qquad\qquad\quad  + 4 n_\kappa' \cdot v_\kappa^{(5)}\Big] dS\,, \\
K_5 =&\ - 4 \int_\Gamma (|\eta_0'|^{-3})' (\nb - b_0)'
(v_\kappa''\cdot \nn)'' dS - 2 \int_\Gamma (|\eta_0'|^{-3})'' (\nb -
b_0) (v_\kappa''\cdot \nn)'' dS\,,\\
K_6 =&\ \int_\Gamma |\eta_0'|^{-3} g_\kappa^{-1} (b_\kappa - b_0)
g_\kappa' (v_\kappa^{(5)} \cdot \nn) dS\,,
\end{align*}
and
\begin{align*}
K_7 =&\ - \int_{\Omega} \nq (\na)_i^j (\zeta_1^2
v_\kappa^{i\prime\prime})_{,j}'' dx\,.
\end{align*}
By (\ref{aetaineq}),
\begin{align*}
|K_1| \le C (1 + t^\delta {\mathcal P}(\nX(T))
\|\nv\|_{H^{2.5}(\Omega)} \|\nv\|_{H^3(\Omega)}
\end{align*}
and hence by Young's inequality and (\ref{regularity1}),
\begin{align}
\Big|\int_0^t K_1 ds\Big| \le&\ C_\epsilon t^\delta \Big[
\|F\|^2_{L^\infty(0,T;L^2(\Omega))} + {\mathcal P}(\nX(T))\Big] +
\epsilon \int_0^t \|\nv\|^2_{H^3(\Omega)} ds \nonumber\\
\le&\ C_\epsilon t^\delta \Big[M_0 + {\mathcal P}(\nX(T))\Big] +
\epsilon \int_0^t \|\nv\|^2_{H^3(\Omega)} ds \,. \label{K1}
\end{align}
Integrating by parts in space, we find that
\begin{align*}
|K_2| \le&\ C \Big[\|\ng - g_0\|_{H^2(\Gamma)}
\|\nv\|_{H^{1.5}(\Gamma)} + \|\ng - g_0\|_{H^{1.5}(\Gamma)}
\|\nv\|_{H^{2.25}(\Gamma)}\Big]\|\neta\|_{H^{3.5}(\Gamma)} \\
& + C \Big[\|\ng - g_0\|_{H^{1.25}(\Gamma)}
\|\neta\|_{H^{3.5}(\Gamma)} + \|\ng - g_0\|_{H^{1.5}(\Gamma)}
\|\neta\|_{H^{3.25}(\Gamma)} \Big]\|\nv\|_{H^{2.5}(\Gamma)}
\end{align*}
and hence
\begin{align}
\Big|\int_0^t K_2 ds\Big| \le C_\epsilon t^\delta {\mathcal
P}(\nX(T)) + \epsilon \int_0^t \Big[\|\nv\|^2_{H^3(\Omega)} +
\|\ng\|^2_{H^{2.5}(\Gamma)}\Big] ds \label{K2}
\end{align}
here we use (\ref{bgsmall}a) to estimate $\|\ng -
g_0\|_{L^\infty(\Gamma)}$. As for $K_3$, integrating by parts for
the first two integrals, we find that
\begin{align*}
& %\Big|\int_\Gamma (\|\eta_0'|^{-3})'(\ng-g_0)' (\eta_\kappa'\cdot
%v_\kappa')'' dS\Big| + \Big|\int_\Gamma (|\eta_0'|^{-3})''(\ng-g_0)
%(\eta_\kappa'\cdot v_\kappa')'' dS\Big| \\
%+&\ \Big| \int_\Gamma \Big[|\eta_0'|^{-3} (g_\kappa - g_0)\Big]''
%g_0'' dS \Big|
|K_3| \le C {\mathcal P}(\nX(T)) \|\ng - g_0\|_{H^2(\Gamma)}
\|\nv\|_{H^2(\Gamma)}\,.
\end{align*}
Therefore,
\begin{align}
\Big|\int_0^t K_3 ds\Big| \le C_\epsilon t^\delta {\mathcal
P}(\nX(T)) + \epsilon \int_0^t \|\nv\|^2_{H^3(\Omega)} ds \,.
\label{K3}
\end{align}
Similarly, integrating by parts and
$H^{1.5}(\Gamma)$-$H^{-1.5}(\Gamma)$ or
$H^{0.5}(\Gamma)$-$H^{-0.5}(\Gamma)$ duality pairing lead to
\begin{align*}
|K_4| + |K_5| \le&\ C {\mathcal P}(\nX(T)) \Big[ \|\nb -
b_0\|_{H^{1.25}(\Gamma)} \|\nn\|_{H^{3.5}(\Gamma)} + \|\nb -
b_0\|_{H^{1.5}(\Gamma)} \|\nn\|_{H^{3.25}(\Gamma)} \\
&\qquad\qquad\quad\ \ + \|\nb - b_0\|_{H^{2.5}(\Gamma)}
\|\nn\|_{H^{2.5}(\Gamma)} \Big] \|\nv\|_{H^{2.5}(\Gamma)}
\end{align*}
and hence
\begin{align}
\Big|\int_0^t (K_4 + K_5)ds\Big| \le C t^\delta {\mathcal P}(\nX(T))
+ \frac{\kappa}{100} \int_0^t \|\eta_\kappa'(s)\|^2_{H^3(\Gamma)}
ds\,, \label{K4K5}
\end{align}
here we use (\ref{bgsmall}a) and (\ref{bkestimate}). By the Leibnitz
rule,
\begin{align*}
K_6 =& - \int_\Gamma \Big[|\eta_0'|^{-3} g_\kappa^{-1} (b_\kappa -
b_0) g_\kappa'\Big]' b_{\kappa t}'' dS - \int_\Gamma
\Big[|\eta_0'|^{-3} g_\kappa^{-1} (b_\kappa - b_0) g_\kappa'\Big]
n_\kappa'\cdot v_\kappa'''' dS \\
& + \int_\Gamma \Big[|\eta_0'|^{-3} g_\kappa^{-1} (b_\kappa - b_0)
g_\kappa'\Big]'\Big[(\eta_\kappa''\cdot n_{\kappa t})'' + v_\kappa
'' \cdot n_\kappa'' + 2 v_\kappa'''\cdot n_\kappa' \Big] dS\,.
\end{align*}
The worst situation for the last integral is when the derivative
outside the bracket is put on $g_\kappa'$. In this case, since
\begin{align}
n_{\kappa t} = - g_\kappa^{-1} (v_\kappa'\cdot\nn) \eta_\kappa'
\quad \text{or}\quad \eta_\kappa'\cdot n_{\kappa t} = -\frac{1}{2}
g_\kappa^{-1} g_\kappa' (v_\kappa'\cdot \nn)\,, \label{etapnt}
\end{align}
the worst term will be
\begin{align*}
\int_\Gamma |\eta_0'|^{-3} g_\kappa^{-1} (b_\kappa - b_0) g_\kappa''
\Big[-\frac{1}{2} g_\kappa^{-1} g_\kappa'''(v_\kappa'\cdot \nn)
-\frac{1}{2} g_\kappa^{-1} g_\kappa' (v_\kappa'''\cdot \nn) + 2
v_\kappa'''\cdot n_\kappa'\Big]dS\,.
\end{align*}
For the first term, since $g_\kappa'' g_\kappa'''$ forms a perfect
derivative, after integrating by parts we have
\begin{align*}
& \Big| \int_\Gamma |\eta_0'|^{-3} g_\kappa^{-1} (b_\kappa - b_0)
g_\kappa'' \Big[-\frac{1}{2} g_\kappa^{-1}
g_\kappa'''(v_\kappa'\cdot \nn)\Big] dS\Big| \\
\le&\ C \mP \|\nb - b_0\|_{H^{1.25}(\Gamma)}
\|\ng\|^2_{H^{2.5}(\Gamma)} \|v_\kappa\|_{H^{1.5}(\Gamma)} \\
& + C \mP \|\nb - b_0\|_{L^\infty(\Gamma)}
\|\ng\|^2_{H^{2.5}(\Gamma)} \|v_\kappa''\cdot \nn\|_{L^2(\Gamma)} \\
& + C \mP \|\nb - b_0\|_{L^\infty(\Gamma)} \|\ng\|_{H^2(\Gamma)}
\|\ng\|_{H^{2.5}(\Gamma)} \|v_\kappa\|_{H^{2.5}(\Gamma)}\,.
\end{align*}
For the remaining two terms, we use
$H^{0.5}(\Gamma)$-$H^{-0.5}(\Gamma)$ duality pairing and obtain
\begin{align*}
& \Big|\int_\Gamma |\eta_0'|^{-3} g_\kappa^{-1} (b_\kappa - b_0)
g_\kappa'' \Big[-\frac{1}{2} g_\kappa^{-1} g_\kappa'
(v_\kappa'''\cdot \nn) + 2 v_\kappa'''\cdot n_\kappa'\Big]dS\Big| \\
\le&\ C \mP \|\nb - b_0\|_{H^1(\Gamma)} \|\ng\|_{H^{2.5}(\Gamma)}
\|\ng\|_{H^{2.5}(\Gamma)} \|\nv\|_{H^{2.5}(\Gamma)}\,.
\end{align*}
Therefore, by (\ref{bgsmall}) and (\ref{bkestimate}), we find that
\begin{align*}
& \Big|\int_0^t \int_\Gamma \Big[|\eta_0'|^{-3} g_\kappa^{-1}
(b_\kappa - b_0) g_\kappa'\Big]'\Big[(\eta_\kappa''\cdot n_{\kappa
t})'' + v_\kappa
'' \cdot n_\kappa'' + 2 v_\kappa'''\cdot n_\kappa' \Big] dS ds\Big| \\
\le&\ C_\epsilon t^\delta P(\nX(T)) + \frac{\kappa}{100} \int_0^t
\|\eta_\kappa'(s)\|^2_{H^3(\Gamma)} ds + \epsilon \int_0^t
\|\nv\|^2_{H^3(\Omega)} ds \,.
\end{align*}
For the second integral of $K_6$, integrating by parts and using
$H^{0.5}(\Gamma)$-$H^{-0.5}(\Gamma)$ duality pairing, we find that
\begin{align*}
& \Big|\int_0^t \int_\Gamma \Big[|\eta_0'|^{-3} g_\kappa^{-1}
(b_\kappa - b_0) g_\kappa'\Big] n_\kappa'\cdot v_\kappa'''' dS
ds\Big| \\
\le&\ C_\epsilon t^\delta P(\nX(T)) + \frac{\kappa}{100} \int_0^t
\|\eta_\kappa'(s)\|^2_{H^3(\Gamma)} ds + \epsilon \int_0^t
\|\nv\|^2_{H^3(\Omega)} ds \,.
\end{align*}
%\begin{align*}
%& \Big|\int_\Gamma |\eta_0'|^{-3} g_\kappa^{-1} (b_\kappa - b_0)
%g_\kappa'(\eta_\kappa''\cdot n_{\kappa t})''' dS \Big|\\
%=&\ \Big|\int_\Gamma \Big[(|\eta_0'|^{-3} g_\kappa^{-1}
%g_\kappa')'(\nb - b_0) + |\eta_0'|^{-3} g_\kappa^{-1} g_\kappa' (\nb
%- b_0)'\Big]\Big[\eta_\kappa'''\cdot n_{\kappa t} +
%\eta_\kappa''\cdot n_{\kappa t}'\Big] dS \Big| \\
%\le&\ C \mP \Big[\|\nb - b_0\|_{L^\infty(\Gamma)}
%\|\ng\|_{H^2(\Gamma)}\|\neta\|_{H^3(\Gamma)}
%\|\nv\|_{H^{2.5}(\Gamma)} \\
%&\qquad\qquad\qquad\quad\ \ + \|\nb - b_0\|_{H^{2.5}(\Gamma)}
%\|\nv\|_{H^{2.5}(\Gamma)} \Big] \\
%\le&\ C {\mathcal P}(\nX(T)) \Big[\|\nb - b_0\|_{L^\infty(\Gamma)}
%\|\neta\|_{H^3(\Gamma)} + \|\nb - b_0\|_{H^{2.5}(\Gamma)}
%\Big]\|\nv\|_{H^{2.5}(\Gamma)}
%\end{align*}
%and hence
%\begin{align*}
%& \Big|\int_0^t \int_\Gamma |\eta_0'|^{-3} g_\kappa^{-1} (b_\kappa -
%b_0) g_\kappa'(\eta_\kappa''\cdot n_{\kappa t})'' dS ds \Big| \\
%\le&\ C t^\delta P(\nX(T)) + \frac{\kappa}{100} \int_0^t
%\|\eta_\kappa'(s)\|^2_{H^3(\Gamma)} ds
%\end{align*}
%Similarly, for the last three terms, we have
%\begin{align*}
%& \Big|\int_0^t \int_\Gamma |\eta_0'|^{-3} g_\kappa^{-1} (b_\kappa -
%b_0) g_\kappa' \Big[v_\kappa''\cdot n_\kappa''' + 3 v_\kappa
%'''\cdot n_\kappa '' + 3 v_\kappa^{(4)}\cdot n_\kappa'\Big] dS
%ds\Big| \\
%\le&\ C t^\delta P(\nX(T)) + \frac{\kappa}{100} \int_0^t
%\|\eta_\kappa'(s)\|^2_{H^3(\Gamma)} ds \,.
%\end{align*}
Finally, for the first integral of $K_6$, we time integrate it first
and then integrate by parts in both space and time to obtain
\begin{align*}
& \int_0^t \int_\Gamma |\eta_0'|^{-3} g_\kappa^{-1} (\nb - b_0)
g_\kappa' b_{\kappa t}''' dS ds \\
=&\ - \int_0^t \int_\Gamma \Big[(|\eta_0'|^{-3} g_\kappa^{-1}
g_\kappa')' (\nb - b_0) + |\eta_0'|^{-3} g_\kappa^{-1} g_\kappa'
(\nb - b_0)'\Big] b_{\kappa t}'' dS ds \\
=&\ - \Big[\int_\Gamma \Big((|\eta_0'|^{-3} g_\kappa^{-1}
g_\kappa')' (\nb - b_0) + |\eta_0'|^{-3} g_\kappa^{-1} g_\kappa'
(\nb - b_0)'\Big) b_\kappa'' dS \Big](t) \\
&\ + \int_0^t \int_\Gamma \Big[(|\eta_0'|^{-3} g_\kappa^{-1}
g_\kappa')' (\nb - b_0) + |\eta_0'|^{-3} g_\kappa^{-1} g_\kappa'
(\nb - b_0)'\Big]_t b_\kappa'' dS ds \,.
\end{align*}
Therefore,
\begin{align*}
& \Big|\int_0^t \int_\Gamma |\eta_0'|^{-3} g_\kappa^{-1} (\nb - b_0)
g_\kappa' b_{\kappa t}''' dS ds\Big| \\
\le&\ \Big|\Big[\int_\Gamma \Big((|\eta_0'|^{-3} g_\kappa^{-1}
g_\kappa')' (\nb - b_0) + |\eta_0'|^{-3} g_\kappa^{-1} g_\kappa'
(\nb - b_0)'\Big) b_\kappa'' dS \Big](t)\Big| \\
&\ + \Big|\int_0^t \int_\Gamma \Big[(|\eta_0'|^{-3} g_\kappa^{-1}
g_\kappa')' (\nb - b_0) + |\eta_0'|^{-3} g_\kappa^{-1} g_\kappa'
(\nb - b_0)'\Big]_t b_\kappa'' dS ds\Big| \\
\le&\ C t^\delta P(\nX(T)) + \frac{\kappa}{100} \int_0^t
\|\eta_\kappa'(s)\|^2_{H^3(\Gamma)} ds\,.
\end{align*}
Combining the estimates above,
\begin{align}
\Big|\int_0^t K_6 ds\Big| \le C t^\delta P(\nX(T)) +
\frac{\kappa}{100} \int_0^t \|\eta_\kappa'(s)\|^2_{H^3(\Gamma)} ds +
\epsilon \int_0^t \|\nv\|^2_{H^3(\Omega)} ds\,. \label{K6}
\end{align}

\begin{remark}\label{need2d}
In the estimate of $K_6$, the fact that $n=2$ is necessary to use the
Sobolev embedding
$L^4(\Gamma)\subset H^{0.25}(\Gamma)$ with
\begin{align}
\|h\|_{H^{0.25}(\Gamma)} \le C \|h\|_{L^4(\Gamma)}
\label{L4Hp25ineq}
\end{align}
and
\begin{align}
\langle fg, h\rangle_\Gamma \le C \|f\|_{H^1(\Gamma)}
\|g\|_{H^{0.5}(\Gamma)} \|h\|_{H^{-0.5}(\Gamma)} \label{Hp5dual}
\end{align}
for some constant $C$. These inequalities no longer holds if $n=3$.
\end{remark}

Now we turn to the estimate of $K_7$. By the identity
\begin{align*}
(\na)_i^j (\zeta_1^2 v_{\kappa}^{i\prime\prime})_{,j} '' &= -
\Big[(a_\kappa'')_i^j (\zeta_1^2 v_{\kappa,j}^i) + 2 (a_\kappa')_i^j
(\zeta_1^2 v_{\kappa,j}^i)' + (\na)_i^j (\zeta_1^2)'' v_{\kappa,j}^i
+ 2 (\na)_i^j (\zeta_1^2)' v_{\kappa,j}^{i\prime}\Big]'' \\
&\ \ \ - (a_\kappa'')_i^j (\zeta_1^2 v_{\kappa,j}^{i\prime\prime}) -
2 (a_\kappa')_i^j (\zeta_1^2)' v_{\kappa,j}^{i\prime} - 2
(a_\kappa')_i^j \zeta_1^2 v_{\kappa,j}^{i\prime\prime} + 2 (\na)_i^j
(\zeta_1 \zeta_{1,j} v_\kappa^{i\prime\prime}) ''\,,
\end{align*}
and inequality (\ref{aetaineq}), we find that
\begin{align*}
& \Big|\int_{\Omega} \nq (\na)_i^j (\zeta_1^2
v_\kappa^{i\prime\prime})_{,j}'' dx\Big| \\
\le&\ C \|\na\|_{H^2(\Omega)}
\Big[\|\nq\|_{H^2(\Omega)}\|\nv\|_{H^{2.5}(\Omega)} +
\|\nq\|_{H^{1.5}(\Omega)} \|\zeta_1
v_{\kappa,j}''\|_{L^2(\Omega)}\Big] \\
\le&\ C_\epsilon (1 + t^\delta {\mathcal P}(\nX(T))
\Big[\|\nv\|^2_{H^2(\Omega)} + \|\nq\|^2_{H^1(\Omega)}\Big] +
\epsilon \Big[\|\nv\|^2_{H^3(\Omega)} +
\|\nq\|^2_{H^2(\Omega)}\Big]\,.
\end{align*}
Therefore, by (\ref{regularity1}) and (\ref{regularity2}),
\begin{align}
\Big|\int_0^t K_7 ds\Big| \le C_\epsilon t^\delta \Big[M_0 +
{\mathcal P}(\nX(T))\Big] + \epsilon \int_0^t
\Big[\|\nv_t\|^2_{H^1(\Omega)} + \|\nv\|^2_{H^3(\Omega)}\Big] ds\,.
\label{K7}
\end{align}
Since
\begin{align*}
\Big|\int_0^t \langle F, (\zeta_1^2 v_\kappa'')'' \rangle ds \Big|
\le C_\epsilon \int_0^t \|F\|^2_{H^1(\Omega)} ds + \epsilon \int_0^t
\|v\|^2_{H^3(\Omega)} ds\,,
\end{align*}
time integrating (\ref{L2H3eq})
together with inequalities (\ref{K1}), (\ref{K2}),
(\ref{K3}), (\ref{K4K5}), (\ref{K6}) and (\ref{K7}), we find that
\begin{align}
&\sup_{0\le t\le T} \Big[\|v_\kappa''\|^2_{L^2(\Omega_1')} + \|\ng -
g_0\|^2_{H^2(\Gamma)} + \|\nb - b_0\|^2_{H^2(\Gamma)} + \kappa \|\neta\|^2_{H^4(\Gamma)} \Big] \label{L2H3ineq}\\
&+ \int_0^T \|v_\kappa''\|^2_{H^1(\Omega_1')} ds \le C_\epsilon
t^\delta \Big[M_0 + {\mathcal P}(\nX(T))\Big] + \epsilon
\int_0^t\Big[\|\nv_t\|^2_{H^1(\Omega)} +
\|\nv\|^2_{H^3(\Omega)}\Big] ds \,. \nonumber
\end{align}
%Since $\|\nv\|_{H^{2.5}(\Gamma)} \le
%C\|v_\kappa''\|_{H^1(\Omega_1)}$ for some constant $C$, by
%(\ref{regularity2}), (\ref{L2H3ineqtemp}) implies
%\begin{align}
%&\sup_{0\le t\le T} \Big[\|v_\kappa''\|^2_{L^2(\Omega_1')} + \|\ng -
%g_0\|^2_{H^2(\Gamma)} + \|\nb - b_0\|^2_{H^2(\Gamma)} + \kappa \|\neta\|^2_{H^4(\Gamma)} \Big] \label{L2H3ineq}\\
%&+ \int_0^T \|\nv\|^2_{H^3(\Omega_0)} ds \le C_\epsilon \Big[M_0 +
%t^\delta{\mathcal P}(\nX(T))\Big] + \epsilon
%\int_0^t\Big[\|\nv_t\|^2_{H^1(\Omega)} +
%\|\nv\|^2_{H^3(\Omega)}\Big] ds \,. \nonumber
%\end{align}

\subsection{Energy estimates for $\nv_t$} Time differentiate
(\ref{weakform}) and then use $\nv_t$ as the test function, we find
that
\begin{align}
& \frac{1}{2}\frac{d}{dt} \Big[\|\nv_t\|^2_{L^2(\Omega)} +
\int_\Gamma |\eta_0'|^{-3} [8|v_\kappa'\cdot \eta_\kappa'|^2 + 2
|v_\kappa''\cdot \nn|^2 ] dS + \kappa \int_\Gamma |\nv^{(4)}|^2 dS\Big] \nonumber\\
+&\ \frac{\nu}{2} \|D_\neta \nv_t\|^2_{L^2(\Omega)} = \langle F_t,
\nv_t \rangle + L_1 + L_2 + L_3 + L_4 + L_5\,, \label{L2H1eq}
\end{align}
where
\begin{align*}
L_1 =&\ - \nu\int_{\Omega} (D_\neta \nv) (a_{\kappa t})_i^k
v_{\kappa t,k}^j dx - \nu \int_{\Omega} \Big[(a_{\kappa t})_i^k
v_{\kappa,k}^j + (a_{\kappa t})_j^k v_{\kappa,k}^i\Big] \na_i^k v_{\kappa t,k}^j dx\,,\\
L_2 =&\ - 2 \int_\Gamma (\ng - g_0) (v_\kappa'\cdot \nv_t') dS + 8
\int_\Gamma |\eta_0'|^{-3} (\eta_\kappa'\cdot v_\kappa')
(v_\kappa'\cdot v_\kappa') dS\,,\\
L_3 =&\ 2 \int_\Gamma |\eta_0'|^{-3} (v_\kappa''\cdot \nn)
(v_\kappa''\cdot n_{\kappa t}) dS - \int_\Gamma |\eta_0'|^{-3}
(\eta_\kappa''\cdot n_{\kappa t}) (v_{\kappa t}'' \cdot \nn) dS \\
&\ - 2\int_\Gamma |\eta_0'|^{-3} (\nb - b_0) (v_{\kappa t}'' \cdot
n_{\kappa t}) dS \,,\\
L_4 =&\ - \int_\Gamma |\eta_0'|^{-3} \Big[g_\kappa^{-1} (\nb-b_0)
g_\kappa' \nn\Big]_t\cdot v_{\kappa t}' dS \,,\\
L_5 =&\ - \int_{\Omega} q_{\kappa t} \na_i^j v_{\kappa t,j}^i dx -
\int_{\Omega} \nq (a_{\kappa t})_i^j v_{\kappa t,j}^i dx\,.
\end{align*}
Since $\na = (\nabla \neta)^{-1}$, $(a_{\kappa t})_i^j = - \na_k^j
\nv_{,\ell}^k \na_i^\ell$ and $\|\na\|_{L^\infty(\Omega)} \le \mP$.
Therefore,
%\begin{align*}
%|L_1| \le C_\epsilon \mP \|\nabla v\|^2_{L^2(\Omega)} + \epsilon
%\|\nv_t\|^2_{H^1(\Omega)}
%\end{align*}
%and hence
\begin{align}
\Big|\int_0^t L_1 ds\Big| \le C_\epsilon t^\delta {\mathcal
P}(\nX(T)) + \epsilon \int_0^t \|\nv_t\|^2_{H^1(\Omega)} ds\,.
\label{L1}
\end{align}
Similar to the estimates in the previous section, integrating by
parts (if necessary) and $H^{1.5}(\Gamma)$-$H^{-1.5}(\Gamma)$ or
$H^{0.5}(\Gamma)$-$H^{-0.5}(\Gamma)$ duality pairing imply that
\begin{align*}
|L_2| \le&\ C \Big[\|\nv\|_{H^{1.5}(\Gamma)} \|v_{\kappa
t}\|_{H^{0.5}(\Gamma)} + \|\nv\|_{H^2(\Gamma)} \|v_{\kappa
t}\|_{L^2(\Gamma)} \Big]\|\ng - g_0\|_{H^{1.25}(\Gamma)} \\
& + C \|\neta\|_{H^{2.5}(\Gamma)}\|\nv\|^2_{H^1(\Gamma)}
\|\nv\|_{H^{2.5}(\Gamma)}
\end{align*}
and hence by (\ref{bgsmall}a),
\begin{align}
\Big|\int_0^t L_2 ds\Big| \le C t^\delta {\mathcal P}(\nX(T))\,.
\label{L2}
\end{align}
For the first integral of $L_3$, we have
\begin{align*}
\Big|\int_\Gamma |\eta_0'|^{-3} (v_\kappa''\cdot \nn)
(v_\kappa''\cdot n_{\kappa t}) dS\Big| \le C_\epsilon
\|\nv\|^2_{H^{2.5}(\Gamma)} \|n_{\kappa t}\|^2_{L^4(\Gamma)} +
\epsilon \|v_\kappa''\cdot \nn\|^2_{L^2(\Gamma)}\,,
\end{align*}
while for the third integral of $L_3$, by integrating by parts,
\begin{align*}
& \int_\Gamma |\eta_0'|^{-3} (\nb - b_0) (v_{\kappa t}'' \cdot
n_{\kappa t}) dS = - \int_\Gamma \Big[|\eta_0'|^{-3} (\nb - b_0)
n_{\kappa t}\Big]'\cdot v_{\kappa t}' dS \\
=& - \int_\Gamma \Big[|\eta_0'|^{-3} (\nb - b_0)\Big]' n_{\kappa t}
\cdot v_{\kappa t}' dS - \int_\Gamma |\eta_0'|^{-3} (\nb - b_0)
n_{\kappa t}'\cdot v_{\kappa t}' dS
\end{align*}
and hence by (\ref{Hp5dual}) and (\ref{bgsmall}a),
\begin{align*}
\Big|\int_\Gamma |\eta_0'|^{-3} (\nb - b_0) (v_{\kappa t}'' \cdot
n_{\kappa t}) dS\Big| \le C_\epsilon {\mathcal P}(\nX(T)) \Big[1+
t^\delta \|v_\kappa''\|^2_{H^1(\Omega_1')}\Big] + \epsilon
\|\nv_t\|^2_{H^1(\Omega)}\,.
\end{align*}
For the second integral of $L_3$,
\begin{align*}
& \int_\Gamma |\eta_0'|^{-3} (\eta_\kappa''\cdot n_{\kappa t})
(v_{\kappa t}'' \cdot \nn) dS = \int_\Gamma |\eta_0'|^{-3}
(\id''\cdot
n_{\kappa t}) (v_{\kappa t}'' \cdot \nn) dS \qquad (\equiv L_{31})\\
&\qquad\qquad\qquad + \int_\Gamma |\eta_0'|^{-3} (\int_0^\cdot
v_\kappa'' ds
\cdot \nn_t) (v_\kappa''\cdot\nn)_t dS \qquad (\equiv L_{32}) \\
&\qquad\qquad\qquad - \int_\Gamma |\eta_0'|^{-3} (\int_0^\cdot
v_\kappa'' ds \cdot \nn_t) (v_\kappa''\cdot\nn_t) dS\,. \qquad
(\equiv L_{33})
\end{align*}
By $H^{1.5}(\Gamma)$-$H^{-1.5}(\Gamma)$ duality pairing,
\begin{align*}
|L_{31}| \le C_\epsilon {\mathcal P}(\nX(T))
\|v_\kappa''\|^2_{H^1(\Omega_1')} + \epsilon
\|\nv_t\|^2_{H^1(\Omega)}
\end{align*}
and standard H$\ddot{\text{o}}$lder's inequality implies
\begin{align*}
|L_{33}| \le C {\mathcal P}(\nX(T)) \|\nv\|_{H^3(\Omega)} \int_0^t
\|\nv\|_{H^3(\Omega)} ds \le C \sqrt{t} {\mathcal P}(\nX(T))
\|\nv\|_{H^3(\Omega)}\,.
\end{align*}
For $L_{32}$, integrating in time and integrating by parts in time
leads to
\begin{align*}
\int_0^t L_{32} ds =& \int_\Gamma |\eta_0'|^{-3} (\int_0^t
v_\kappa'' ds \cdot \nn_t) (v_\kappa''\cdot\nn) dS - \int_0^t
\int_\Gamma |\eta_0'|^{-3} (v_\kappa'' \cdot \nn_t) (v_\kappa''\cdot\nn) dS ds\\
&- \int_0^t \int_\Gamma |\eta_0'|^{-3} (\int_0^\cdot v_\kappa'' ds
\cdot \nn_{tt}) (v_\kappa''\cdot\nn) dS ds\,.
\end{align*}
The first two integrals can be bounded by
\begin{align*}
C \sqrt{t} {\mathcal P}(\nX(T)) + C {\mathcal P}(\nX(T)) \int_0^t
\|v_\kappa''\|^2_{H^3(\Omega_1')} ds\,.
\end{align*}
The worst term in the last integral is
\begin{align*}
\int_0^t \int_\Gamma |\eta_0'|^{-3} (\int_0^\cdot v_\kappa'' ds
\cdot \eta_\kappa') g_\kappa^{-1} (v_{\kappa t}'\cdot \nn)
(v_\kappa''\cdot\nn) dS ds\,,
\end{align*}
and integrating by parts, the worst term becomes
\begin{align}
& \int_0^t \int_\Gamma |\eta_0'|^{-3} (\int_0^\cdot v_\kappa''' ds
\cdot \eta_\kappa') g_\kappa^{-1} (v_{\kappa t}\cdot \nn)
(v_\kappa''\cdot\nn) dS ds \nonumber\\
=& \int_0^t \int_\Gamma |\eta_0'|^{-3} [(\eta_\kappa''' - \id''')
\cdot \eta_\kappa'] g_\kappa^{-1} (v_{\kappa t}\cdot \nn)
(v_\kappa''\cdot\nn) dS ds\,. \label{goodtransform}
\end{align}
Standard H$\ddot{\text{o}}$lder's inequality and interpolation
inequalities show that
\begin{align*}
& \Big|\int_0^t \int_\Gamma |\eta_0'|^{-3} [(\eta_\kappa''' -
\id''') \cdot \eta_\kappa'] g_\kappa^{-1} (v_{\kappa t}\cdot \nn)
(v_\kappa''\cdot\nn) dS ds \Big| \\
\le&\ {\mathcal P}(\nX(T)) \int_0^t (1 +
\|\neta\|_{H^{3.25}(\Gamma)}) \|\nv_t\|_{H^{0.25}(\Gamma)} ds \\
\le&\ C_\epsilon t^\delta {\mathcal P}(\nX(T)) + \epsilon \int_0^t
\|\nv_t\|^2_{H^1(\Omega)}\,.
\end{align*}
Estimates for $L_{31}$, $L_{32}$ and $L_{33}$ leads to
\begin{align*}
& \Big| \int_0^t \int_\Gamma |\eta_0'|^{-3} (\eta_\kappa''\cdot
n_{\kappa t}) (v_{\kappa t}'' \cdot \nn) dS ds\Big| \\
\le&\ C_\epsilon t^\delta {\mathcal P}(\nX(T)) + C_\epsilon
{\mathcal P}(\nX(T)) \int_0^t \|v_\kappa''\|^2_{H^1(\Omega_1')} ds +
\epsilon \int_0^t \|\nv_t\|^2_{H^1(\Omega)} ds \,.
\end{align*}
Therefore,
\begin{align}
\Big|\int_0^t L_3 ds\Big| \le&\ C_\epsilon t^\delta {\mathcal
P}(\nX(T)) + C_\epsilon {\mathcal P}(\nX(T)) \int_0^t
\|v_{\kappa}''\|^2_{H^1(\Omega_1')} ds \label{L3} \\
& + \epsilon \int_0^t \|v_\kappa''\cdot \nn\|^2_{L^2(\Gamma)} ds +
\epsilon \int_0^t \|\nv_t\|^2_{H^1(\Omega)} ds + \frac{\kappa}{100}
\int_0^t \|\eta_{\kappa}'\|^2_{H^3(\Gamma)} ds \,. \nonumber
\end{align}

$L_4$ can be rewritten as
\begin{align*}
L_4 = - \int_\Gamma |\eta_0'|^{-3} (\nb - b_0) (g_\kappa^{-1}
g_\kappa' \nn)_t\cdot v_{\kappa t}' dS - \int_\Gamma |\eta_0'|^{-3}
(\nb - b_0)_t g_\kappa^{-1} g_\kappa' \nn \cdot v_{\kappa t}' dS\,.
\end{align*}
The presence of $(\nb - b_0)$ and
$H^{0.5}(\Gamma)$-$H^{-0.5}(\Gamma)$ duality pairing imply that
\begin{align}
\Big|\int_0^t \int_\Gamma |\eta_0'|^{-3} (\nb - b_0) (g_\kappa^{-1}
g_\kappa' \nn)_t\cdot v_{\kappa t}' dS ds\Big| \le&\ C t^\delta
{\mathcal P}(\nX(T))\,. \label{L41}
\end{align}
For the second integral, since $b_{\kappa t} = v_\kappa''\cdot \nn +
\eta_\kappa''\cdot \nn_t$, by $H^{0.5}(\Gamma)$-$H^{-0.5}(\Gamma)$
duality pairing,
\begin{align}
& \Big|\int_0^t \int_\Gamma |\eta_0'|^{-3} (\nb - b_0)_t
g_\kappa^{-1} g_\kappa' (\nn \cdot v_{\kappa t}') dS ds\Big| \nonumber\\
\le&\ C \int_0^t {\mathcal P}(\|\eta\|^2_{H^{2.5}(\Gamma)}) (1 +
\|\ng\|_{H^2(\Gamma)}) \|v_\kappa''\|_{H^1(\Omega_1')}
\|\nv_t\|_{H^1(\Omega)} ds \nonumber\\
\le&\ C_\epsilon {\mathcal P}(\nX(T)) \int_0^t
\|v''\|^2_{H^1(\Omega_1')} ds + \epsilon \int_0^t
\|\nv_t\|^2_{H^1(\Omega)} ds\,. \label{L42}
\end{align}
Combining (\ref{L41}) and (\ref{L42}), we have
\begin{align}
\Big|\int_0^t L_4 ds\Big| \le&\ C t^\delta {\mathcal P}(\nX(T)) +
C_\epsilon {\mathcal P}(\nX(T)) \int_0^t \|v''\|^2_{H^1(\Omega_1')}
ds \label{L4} \\
& + \epsilon \int_0^t \|\nv_t\|^2_{H^1(\Omega)} ds\,. \nonumber
\end{align}
For $L_5$, first note that by $(a_{\kappa t})_i^j = - (\na)_k^j
\nv_{,\ell}^k (\na)_i^\ell$, we find that
\begin{align}
& \Big|\int_{\Omega} \nq (a_{\kappa t})_i^j v_{\kappa t,j}^i dx\Big|
\le \|\nq\|_{H^{0.5}(\Omega)} \|\na\|^2_{L^\infty(\Omega)}
\|\nv\|_{H^{1.5}(\Omega)} \|\nv_t\|_{H^1(\Omega)} \nonumber\\
\le&\ C_\epsilon \mP \|\nq\|^2_{L^2(\Omega)} + \epsilon
\Big[\|\nq\|^2_{H^1(\Omega)} + \|\nv_t\|^2_{H^1(\Omega)}\Big]
\label{L51}\,.
\end{align}
By the ``divergence free'' condition, $(\na)_i^j v_{\kappa t,j}^i =
- (a_{\kappa t})_i^j v_{\kappa,j}^i$. Time integrating the first
integral,
\begin{align*}
& \int_0^t \int_{\Omega} q_{\kappa t} (\na)_i^j v_{\kappa t,j}^i dx
ds \\
=& - \int_0^t \int_{\Omega} \nq \Big[(\na)_k^j \nv_{,\ell}^k
(\na)_i^\ell v_{\kappa,j}^i\Big]_t dx ds + \int_{\Omega} \nq(s)
(a_{\kappa t}(s))_i^j v_{\kappa,j}^i(s) dx\Big|_{s=0}^{s=t}\,.
\end{align*}
Similar to (\ref{L51}), we have
\begin{align}
& \Big|\int_0^t \int_{\Omega} \nq \Big[(\na)_k^j \nv_{,\ell}^k
(\na)_i^\ell v_{\kappa,j}^i\Big]_t dx ds\Big| \le C_\epsilon \mP
\int_0^t \|\nq\|^2_{L^2(\Omega)}ds \label{L52}\\
&\qquad\qquad\qquad\qquad\qquad + \epsilon \int_0^t
\Big[\|\nq\|^2_{H^1(\Omega)} + \|\nv_t\|^2_{H^1(\Omega)}\Big] ds \,.
\nonumber
\end{align}
For the boundary term, when $s=0$, it is bounded by a constant
independent of $\kappa$, say $M_0$. When $s=t$, first note that by
$(a_{\kappa t})_i^j = - \na_k^j \nv_{,\ell}^k \na_i^\ell$,
\begin{align*}
& \|(a_{\kappa t})_i^j (t) - (a_{\kappa t})_i^j
(0)\|^2_{L^4(\Omega)} \\
\le&\ \|(\na_k^j - \delta_k^j) v_{\kappa,\ell}^k
\na_i^\ell\|^2_{L^4(\Omega)} + \|(v_{\kappa,\ell}^j - u_{0,\ell}^j)
\na_i^\ell\|^2_{L^4(\Omega)} + \|u_{0,\ell}^j (\na_i^\ell -
\delta_i^\ell)\|^2_{L^4(\Omega)} \\
\le&\ C (t + \sqrt{t}) {\mathcal P}(\nX(T))\,.
\end{align*}
Therefore, by adding and subtracting $\displaystyle{\int_{\Omega}
\nq(t) (a_{\kappa t}(0))_i^j v_{\kappa,j}^i(t) dx}$, we find that
\begin{align}
& \Big|\int_{\Omega} \nq(t) (a_{\kappa t}(t))_i^j v_{\kappa,j}^i(t)
dx\Big| \nonumber\\
\le&\ \Big|\int_{\Omega} \nq(t) (a_{\kappa t}(t) - a_{\kappa
t}(0))_i^j v_{\kappa,j}^i(t) dx\Big| + \Big|\int_{\Omega} \nq(t)
u_{0,i}^j v_{\kappa,j}^i dx\Big| \nonumber\\
\le&\ C t^\delta {\mathcal P}(\nX(T)) + C_\epsilon \Big[
\|u_0\|^2_{H^1(\Omega)} + \sqrt{t} {\mathcal P}(\nX(T))\Big] +
\epsilon \|q\|^2_{L^2(\Omega)}\,, \label{L53}
\end{align}
where (\ref{H1estimate1}) is used to estimate
$\|v\|^2_{H^1(\Omega)}$. Combining (\ref{L51}), (\ref{L52}) and
(\ref{L53}), by (\ref{qL2estimate1}) and (\ref{regularity1}) we
obtain
\begin{align}
\Big|\int_0^t L_5 ds\Big| \le&\ C_\epsilon t^\delta {\mathcal
P}(\nX(T)) + \epsilon \int_0^t \|\nv_t\|^2_{H^1(\Omega)} ds +
\epsilon \|\nv_t\|^2_{L^2(\Omega)} \,.\label{L5}
\end{align}
Time integrating (\ref{L2H1eq}), choosing $\epsilon>0$ small enough
together with inequalities (\ref{L1}), (\ref{L2}), (\ref{L3}),
(\ref{L4}) and (\ref{L5}), we find that
\begin{align}
& \sup_{0\le t\le T}\Big[\|\nv_t\|^2_{L^2(\Omega)} +
\|v_\kappa'\cdot \eta_\kappa'\|^2_{L^2(\Gamma)} + \|v_\kappa''\cdot
\nn\|^2_{L^2(\Gamma)} + \kappa \|\nv\|^2_{H^4(\Gamma)} \Big] \label{L2H1ineq} \\
& + \int_0^T \|v_t\|^2_{H^1(\Omega)} ds \le M_0 + C T^\delta
{\mathcal P}(\nX(T)) + C {\mathcal P}(\nX(T)) \int_0^t
\|v_\kappa''\|^2_{H^1(\Omega_1')} ds \nonumber\\
&\qquad\qquad\qquad\qquad\quad + \frac{\kappa}{100} \int_0^T
\|\eta_{\kappa}'\|^2_{H^3(\Gamma)} ds \,, \nonumber
\end{align}
for some constant $M_0$ depending on $\|u_0\|^2_{H^2(\Omega)}$,
$\|\id\|^2_{H^2(\Gamma)}$, $\|F\|^2_{L^\infty(0,T;H^1(\Omega))}$ and
$\|F_t\|^2_{L^2(0,T;L^2(\Omega))}$.

\subsection{$\kappa$-independent estimates} Let $\bbE_1$ and $\bbE_2$ be the left-hand side quantities of
(\ref{L2H3ineq}) and (\ref{L2H1ineq}), respectively. Then
\begin{subequations}\label{E1E2}
\begin{align}
\bbE_1 &\le C_\epsilon T^\delta \Big[M_0 + {\mathcal P}(\nX(T))\Big]
+ \epsilon \bbE_2 + \epsilon \int_0^t \|\nv\|^2_{H^3(\Omega)} ds \,,\\
\bbE_2 &\le M_0 + C T^\delta {\mathcal P}(\nX(T)) + C {\mathcal
P}(\nX(T)) \bbE_1 + \frac{\kappa}{100} \bbE_1\,.
\end{align}
\end{subequations}
By (\ref{E1E2}a), for $\epsilon >0$ small enough (but not fixed yet,
say $C\epsilon {\mathcal P}(\nX(T)) \le 0.5$), we have
\begin{align}
\bbE_2 &\le M_0 + C T^\delta {\mathcal P}(\nX(T)) + \epsilon
\int_0^t \|\nv\|^2_{H^3(\Omega)} + \frac{\kappa}{100} \bbE_1\,.
\label{bbE2}
\end{align}
Finally, since
\begin{align}
\|\nv\|^2_{H^2(\Omega)} \le&\ C_\epsilon \|\nv\|^2_{H^1(\Omega)} +
\epsilon \|\nv\|^2_{H^3(\Omega)} \nonumber\\
\le&\ C_\epsilon \Big[\|u_0\|^2_{H^1(\Omega)} + \bbE_2 \Big] +
\epsilon \int_0^t \|v\|^2_{H^3(\Omega)} ds \label{H2estimate}
\end{align}
combining (\ref{bkestimate}), (\ref{gkestimate}),
(\ref{regularity2}), (\ref{E1E2}a), (\ref{bbE2}) and
(\ref{H2estimate}), and choosing $\epsilon>0$ small enough, we have
\begin{align}
%&\sup_{0\le t\le T} \Big[\|\nv_t\|^2_{L^2(\Omega)} +
%\|\nv\|^2_{H^2(\Omega)} + \|\ng\|^2_{H^2(\Gamma)} + \|\nb\|^2_{H^2(\Gamma)}
%+ \|v_\kappa'\cdot \eta_\kappa'\|^2_{L^2(\Gamma)} \nonumber\\
%&\qquad\quad + \|v_\kappa''\cdot \nn\|^2_{L^2(\Gamma)} + \kappa
%\|\neta\|^2_{H^4(\Gamma)} + \kappa \|\nv\|^2_{H^4(\Gamma)} \Big] + \|\nv\|^2_{{\mathcal V}^3(T)} \\
%& + \|\nb\|^2_{L^2(0,T;H^{2.5}(\Gamma))} +
%\|\ng\|^2_{L^2(0,T;H^{2.5}(\Gamma))}
\nX(T) \le C \Big[M_0 + T^\delta{\mathcal P}(\nX(T))\Big]\,.
\label{mainestimate}
\end{align}
%Define
%\begin{align*}
%Y_\kappa(T) =& \sup_{0\le t\le T} \Big[\|\nv_t\|^2_{L^2(\Omega)} +
%\|\nv\|^2_{H^2(\Omega)} + \|\ng\|^2_{H^2(\Gamma)} +
%\|\nb\|^2_{H^2(\Gamma)}
%+ \|v_\kappa'\cdot \eta_\kappa'\|^2_{L^2(\Gamma)} \nonumber\\
%&\qquad\quad + \|v_\kappa''\cdot \nn\|^2_{L^2(\Gamma)} + \kappa
%\|\neta\|^2_{H^4(\Gamma)} + \kappa \|\nv\|^2_{H^4(\Gamma)}\Big]
%+ \|\nv\|^2_{{\mathcal V}^3(T)} \\
%& + \|\nb\|^2_{L^2(0,T;H^{2.5}(\Gamma))} +
%\|\ng\|^2_{L^2(0,T;H^{2.5}(\Gamma))} \,.
%\end{align*}
$X_\kappa$ is clearly continuous in its variable. By Lemma
\ref{polytypeineq}, there is a constant $M$ independent of $\kappa$
and $T_1\le T_\kappa$ so that
\begin{align*}
X_\kappa(t) \le M \qquad\forall\ t\in [0,T_1]\,.
\end{align*}
Without loss of generality, we may assume that $T_1 = T_\kappa$ (by
setting $T_\kappa$ equaling $T_1$). Let $X_{0.5}(t) \le M$ for $t\in
[0,T_{0.5}]$. For $\kappa < 0.5$, say, $\kappa = 0.1$, $X_{0.1}(t)
\le M$ for $t\in [0,T_{0.1}]$ where $T_{0.1}$ is in general smaller
than $T_{0.5}$. Since this estimate is independent of $\kappa$, we
are able to extend the time interval $[0,T_{0.1}]$ in which the
fixed point $v_{0.1}$, $\eta_{0.1}$, $g_{0.1}$ and $b_{0.1}$ exist.
This extension will proceed until $T_{0.1}$ hits $T_{0.5}$, and
hence $X_{0.1}(t) \le M$ for $t\in[0,T_{0.5}]$. This argument holds
for all $\kappa < 0.5$, so we conclude that (with $T \equiv
T_{0.5}$)
\begin{align}
X_\kappa(t) \le M \qquad\forall\ t\in [0,T],\quad \kappa\in
(0,0.5]\,. \label{kappaindependentestimate}
\end{align}

\begin{remark}
By (\ref{regularity2}), we can also include
$\|\nq\|^2_{L^2(0,T;H^2(\Omega))}$ in $\nX(T)$.
\end{remark}

\subsection{Weak limits of $\nv$ as $\kappa\to 0$} By
(\ref{kappaindependentestimate}), there exist $v$ (and $v_t$) so
that
\begin{subequations}\label{weakconvergence1}
\begin{align}
v_{\kappa_i} &\rightharpoonup v \qquad\ \text{in } L^2(0,T;H^3(\Omega))\,, \\
v_{\kappa_i} &\to v \qquad\ \text{in } L^2(0,T;H^2(\Omega))\,, \\
v_{\kappa_i t} &\rightharpoonup v_t \qquad \text{in } L^2(0,T;H^1(\Omega))\,,\\
v_{\kappa_i t} &\to v_t \qquad \text{in } L^2(0,T;L^2(\Omega))\,,
\end{align}
\end{subequations}
for some subsequence $v_{\kappa_i}$. Also, there exists $\eta$ (the
associated Lagrangian variable of $v$), $g$, $b$ and $n$ so that
\begin{subequations}\label{weakconvergence2}
\begin{align}
\eta_{\kappa_i} &\to \eta \qquad \ \text{in } L^2(0,T;H^3(\Gamma))\cap L^2(0,T;H^2(\Omega))\,,\\
g_{\kappa_i} &\to g \qquad\ \text{in } L^2(0,T;H^2(\Gamma))\,,\\
b_{\kappa_i} &\to b \qquad\ \text{in } L^2(0,T;H^2(\Gamma))\,,\\
n_{\kappa_i} &\to n \qquad\ \text{in } L^2(0,T;H^2(\Gamma))\,.
\end{align}
\end{subequations}
Since $\eta_{\kappa_i}$ converges a.e. to $\eta$ in $H^3(\Gamma)$,
we have that $g = |\eta'|^2$, $b = \eta''\cdot n$, and
$n = (-\eta_2', \eta_1')/|\eta'|$. Also, since $a_{\kappa_i} \to a$
strongly in $L^2(0,T;H^1(\Omega))$, by (\ref{weakform}) we conclude
that $v$, $q$, $\eta$, $g$, $b$ satisfy
(\ref{solidshelleq}).

\section{Uniqueness} Let $v$ and $\tv$ in ${\mathcal V}^5(T)$ be two solutions to
(\ref{solidshelleq}) ($q$ and $\tq\in L^2(0,T;H^4(\Omega))\cap
L^\infty(0,T;H^3(\Omega))$, $q_t$ and $\tq_t\in
L^2(0,T;H^2(\Omega))$, $g$, $\tg$, $b$, $\tb \in
L^2(0,T;H^{4.5}(\Gamma)) \cap L^\infty(0,T;H^4(\Gamma))$
), and $w = v-\tv$, $r=q-\tq$, $\bbE = \eta - \teta$. %they satisfy
%\begin{align*}
%& \langle v_t,\phi\rangle + \frac{\nu}{2} \langle a_i^k v_{,k}^j +
%a_j^k v^j_{,k}, a_i^k \phi_{,k}^i +
%a_j^k \phi_{,k}^j\rangle + \langle q, a_i^j\phi_{,j}^i\rangle \\
%+& \int_\Gamma |\eta_0'|^{-3}\Big[ 4 (g - g_0) (\eta'\cdot \phi') + 2 (b - b_0) (\phi''\cdot n)\Big] dS \\
%-& \int_\Gamma |\eta_0'|^{-3} g^{-1} (b - b_0) g' (\phi'\cdot n) dS
%= \langle f\circ\eta,\phi\rangle
%\end{align*}
%and
%\begin{align*}
%& \langle \tv_t,\phi\rangle + \frac{\nu}{2} \langle \ta_i^k
%\tv_{,k}^j + \ta_j^k \tv^j_{,k}, \ta_i^k \phi_{,k}^i +
%\ta_j^k \phi_{,k}^j\rangle + \langle \tq, \ta_i^j\phi_{,j}^i\rangle \\
%+& \int_\Gamma |\eta_0'|^{-3}\Big[ 4 (\tg - g_0) (\eta'\cdot \phi') + 2 (\tb - b_0) (\phi''\cdot \tn)\Big] dS \\
%-& \int_\Gamma |\eta_0'|^{-3} \tg^{-1} (\tb - b_0) \tg' (\phi'\cdot
%\tn) dS = \langle f\circ\teta,\phi\rangle\,,
%\end{align*}
%where all the $\tilde{\ }$-variables on the boundary are constructed
%from the associated Lagrangian variable $\teta$ of $\tv$. Let $w= v
%- \tv$, then the difference of these two equalities is
Then $w$, $r$, $\bbE$ satisfy
\begin{subequations}\label{diffeq}
\begin{alignat}{2}
w^i_t - \nu (a_\ell^j a_\ell^k w_{,k}^i)_{,j} =& - a_i^j r_{,j} + (\delta F)^i && \text{in }(0,T)\times\Omega\,,\\
a_i^j w_{,j}^i =&\ \delta a && \text{in }(0,T)\times\Omega \,,\\
\Big[\nu (a_i^k w_{,k}^j + a_j^k w^i_{,k}) - r
\delta_{ij}\Big]a_j^\ell N_\ell =&\  \tL(\bbE) + \delta L_1 + \delta L_2 \quad&& \text{on }(0,T)\times\Gamma\,, \\
& + \delta L_3 + \delta L_4 && \nonumber\\
w(0,x)=&\ 0 && \text{in }\Omega\,,
\end{alignat}
\end{subequations}
where
\begin{align*}
\delta F =&\ f\circ\eta - f\circ\teta + \nu [(a_\ell^k a_\ell^j -
\ta_\ell^k\ta_\ell^j) \tv_{,j}^i]_{,k} + \nu [(a_\ell^k a_i^j -
\ta_\ell^k\ta_i^j) \tv_{,j}^\ell]_{,k} - (a_i^k -
\ta_i^k)\tq_{,k}\,, \\
\tL(\bbE) =&\ - 4\Big[|\eta_0'|^{-3} (\bbE\cdot \eta')\eta'\Big]' + 2 \Big[|\eta_0'|^{-3} (\bbE''\cdot n) n\Big]''\,,\\
\delta a =&\ - (a_i^j - \ta_i^j)\tv_{,j}^i\,,\\
\delta L_1 =&\ -\nu \Big[(a_i^k - \ta_i^k)\tv_{,k}^j + (a_j^k -
\ta_j^k) \tv_{,k}^i \Big] a_j^\ell N_\ell - \nu \Big[\ta_i^k
\tv_{,k}^j + \ta_j^k \tv_{,k}^i \Big] (a_j^\ell - \ta_j^\ell) N_\ell
\\
&\ - \tq (a_i^\ell - \ta_i^\ell) N_\ell\,, \\
\delta L_2 =&\ 4 \Big[|\eta_0'|^{-3}[(\teta'\cdot\bbE')\eta' +
(\tg-g_0) \bbE']\Big]'\,, \\
\delta L_3 =&\ 2 \Big[|\eta_0'|^{-3} [ \teta''\cdot (n-\tn) n + (\tb
- b_0) (n-\tn) ] \Big]''\,, \\
\delta L_4 =&\ \Big[|\eta_0'|^{-3}[(g^{-1} - \tg^{-1})(b-b_0) g'n
+ \tg^{-1} (b-\tb) g' n + \tg (\tb - b_0) (g - \tg)' n \\
&\qquad\quad\ + \tg^{-1}(\tb - b_0) \tg' (n-\tn)]\Big]'\,,
\end{align*}
with the following inequalities from \cite{ChCoSh2006}:
\begin{subequations}\label{deltaF}
\begin{align}
%\|a-\ta\|^2_{H^{k+1}(\Omega)} \le&\ C t \int_0^t
%\|w\|^2_{H^{k+2}(\Omega)} ds\,, \quad\text{for $k=0$, $1$,}\\
\|\delta a\|^2_{H^k(\Omega)} + \|\delta F\|^2_{H^k(\Omega)} \le&\
Ct\int_0^t \|w\|^2_{H^{k+1}(\Omega)} ds \,,\quad\text{for $k=0$, $1$, $2$,}\\
 \|(\delta a)_t\|^2_{L^2(\Omega)} + \|(\delta F)_t\|^2_{L^2(\Omega)} %\le&\ C\Big[t
%(1+\|q_t\|^2_{H^1(\Omega)} )\int_0^t \|w\|^2_{H^3(\Omega)} ds +
%\|w\|^2_{H^2(\Omega)}\Big]\,, \\
\le&\ C \sqrt{t}\int_0^t \|w\|^2_{{\mathcal V}^3(\Omega)} ds\,,\\
\|\delta L_1\|^2_{H^{1.5}(\Gamma)} + \|(\delta
L_1)_t\|^2_{L^2(\Gamma)} \le&\ C t \int_0^t \|w\|^2_{H^3(\Omega)}
ds\,.
\end{align}
\end{subequations}
Furthermore, $w$, $r$ and $\bbE$ satisfy the following variational
form:
\begin{align}
& \langle w_t,\phi\rangle + \frac{\nu}{2} \langle a_i^k w_{,k}^j +
a_j^k w^j_{,k}, a_i^k \phi_{,k}^i +
a_j^k \phi_{,k}^j\rangle + \langle r, a_i^j\phi_{,j}^i\rangle + \int_\Gamma (\delta L_1)\cdot \phi dS \nonumber\\
+& \int_\Gamma |\eta_0'|^{-3}\Big[ 4 (g - \tg) (\eta'\cdot \phi') + 2 (b - \tb) (\phi''\cdot n)\Big] dS \nonumber\\
+& \int_\Gamma |\eta_0'|^{-3}\Big[ 4 (\tg - g_0) (\bbE'\cdot \phi') + 2 (\tb - b_0) (\phi''\cdot (n-\tn))\Big] dS \label{weq}\\
-& \int_\Gamma |\eta_0'|^{-3} (g^{-1} - \tg^{-1}) (b - b_0) g'
(\phi'\cdot n) dS - \int_\Gamma |\eta_0'|^{-3} \tg^{-1} (b - \tb) g'
(\phi'\cdot n) dS \nonumber\\
-&\ \int_\Gamma |\eta_0'|^{-3} \tg^{-1} (\tb - b_0) (g-\tg)'
(\phi'\cdot n) dS - \int_\Gamma |\eta_0'|^{-3} \tg^{-1} (\tb - b_0)
\tg' (\phi'\cdot (n-\tn)) dS \nonumber\\
=&\ \langle \delta F, \phi\rangle \nonumber
\end{align}
for all $\phi\in\H1H2$%, where all the $\tilde{\ }$-variables on the
%boundary are constructed from the associated Lagrangian variable
%$\teta$ of $\tv$
.

\subsection{Some a priori estimates} Similar to
(\ref{regularity1temp}) and (\ref{regularity2temp}), solving a
Stokes problem (formed from (\ref{diffeq}a) and (\ref{diffeq}b))
gives us
\begin{align}
\|w\|^2_{H^2(\Omega)} + \|r\|^2_{H^1(\Omega)} \le&\ C \Big[\|\delta
F\|^2_{L^2(\Omega)} + \|w_t\|^2_{L^2(\Omega)} + \|\delta
a\|^2_{H^1(\Omega)} + \|w\|^2_{H^{1.5}(\Gamma)}\Big] \nonumber \\
\le&\ C \Big[\|w_t\|^2_{L^2(\Omega)} + \|w\|^2_{H^{1.5}(\Gamma)} +
t\int_0^t \|w\|^2_{H^2(\Omega)} ds\Big] \label{wrregularity1}
\end{align}
and
\begin{align}
\|w\|^2_{H^3(\Omega)} + \|r\|^2_{H^2(\Omega)} \le C
\Big[\|w_t\|^2_{H^1(\Omega)} + \|w\|^2_{H^{2.5}(\Gamma)} + t\int_0^t
\|w\|^2_{H^3(\Omega)} ds\Big]\,. \label{wrregularity2}
\end{align}
For $T$ small enough, (\ref{wrregularity2}) implies that
\begin{align}
\int_0^t \Big[\|w\|^2_{H^3(\Omega)} + \|r\|^2_{H^2(\Omega)}\Big] ds
\le C \int_0^t \Big[\|w_t\|^2_{H^1(\Omega)} +
\|w\|^2_{H^{2.5}(\Gamma)}\Big] ds \,. \label{wrregularity3}
\end{align}

We can also setup elliptic equations for $b - \tb$ and $g - \tg$ and
obtain the following elliptic estimates for $b-\tb$ and $g - \tg$
(where we use (\ref{wrregularity1}) and (\ref{wrregularity2}) to
estimate the norm of $r$):
\begin{subequations}\label{bgdiffineq}
\begin{align}
\|b-\tb\|^2_{H^{1.5}(\Gamma)} %\le&\ C t \Big[\int_0^t
%\|w\|^2_{H^3(\Omega)} ds + \|w\|^2_{H^2(\Omega)} \Big] + C
%\Big[\|r\|^2_{H^1(\Omega)} +
%\|w\|^2_{H^2(\Omega)}\Big] \nonumber\\
\le&\ C t \int_0^t
\|w\|^2_{H^3(\Omega)} ds + C \Big[\|w_t\|^2_{L^2(\Omega)} + \|w\|^2_{H^{1.5}(\Gamma)}\Big] \,, \\
\|b-\tb\|^2_{H^{2.5}(\Gamma)} %\le&\ C t \Big[\int_0^t
%\|w\|^2_{H^3(\Omega)} ds + \|w\|^2_{H^3(\Omega)} +
%\|g - \tg\|^2_{H^{2.5}(\Gamma)} \Big] \nonumber\\
%& + C \Big[\|r\|^2_{H^1(\Omega)} + \|w\|^2_{H^2(\Omega)}\Big]
%\nonumber\\
\le&\ C t \Big[\int_0^t \|w\|^2_{H^3(\Omega)} ds +
\|w\|^2_{H^3(\Omega)} + \|g - \tg\|^2_{H^{2.5}(\Gamma)} \Big] \\
& + C \Big[\|w_t\|^2_{L^2(\Omega)} + \|w\|^2_{H^{1.5}(\Gamma)} \Big] \,, \nonumber\\
\|g-\tg\|^2_{H^{2.5}(\Gamma)} \le&\ C \Big[ \|w\|^2_{H^3(\Omega)} +
\|b - \tb\|^2_{H^{2.5}(\Gamma)} + \|\tb -
b_0\|^2_{H^{2.5}(\Gamma)} \|b - \tb\|^2_{H^{1.5}(\Gamma)} \Big] \nonumber\\
& + C t \int_0^t \|w\|^2_{H^3(\Omega)} ds\,.
\end{align}
\end{subequations}
Using (\ref{bgdiffineq}b) in (\ref{bgdiffineq}c), for $T$ small
enough, $\|g-\tg\|^2_{H^{2.5}(\Gamma)}$ dependence on the right-hand
side can be absorbed by the left-hand side of (\ref{bgdiffineq}c),
so we conclude that
\begin{align*}
\int_0^t \|g-\tg\|^2_{H^{2.5}(\Gamma)} ds &\le C \int_0^t
\Big[\|w_t\|^2_{H^1(\Omega)} + \|w\|^2_{H^{2.5}(\Gamma)} \Big] ds\,, \\
\int_0^t \|b-\tb\|^2_{H^{2.5}(\Gamma)} ds &\le C \int_0^t
\Big[t\|w_t\|^2_{H^1(\Omega)} + t\|w\|^2_{H^{2.5}(\Gamma)} +
\|w_t\|^2_{L^2(\Omega)} + \|w\|^2_{H^{1.5}(\Gamma)} \Big] ds \,.
\end{align*}
By the identity $\displaystyle{\eta'' = g^{-1} (\eta''\cdot \eta')
\eta' + (\eta''\cdot n) n = \frac{1}{2} g^{-1} g' \eta' + b n}$, we
find that
\begin{align*}
\|\bbE\|^2_{H^{3.5}(\Gamma)} \le C \Big[ \|g -
\tg\|^2_{H^{2.5}(\Gamma)} + \|\bbE\|^2_{H^{2.5}(\Gamma)} + \|b -
\tb\|^2_{H^{1.5}(\Gamma)}\Big]
\end{align*}
and hence
\begin{align}
\int_0^t \|\bbE\|^2_{H^{3.5}(\Gamma)} ds \le C \int_0^t\Big[
\|w\|^2_{H^{2.5}(\Gamma)} + \|w_t\|^2_{H^1(\Omega)} \Big]ds\,.
\label{bbEH3p5estimate}
\end{align}
By (\ref{normalderivatives}b),
\begin{align*}
\|n - \tn\|^2_{H^s(\Gamma)} \le C \Big[\|\bbE\|^2_{H^s(\Gamma)} +
\|b-\tb\|^2_{H^{s-1}(\Gamma)}\Big] \quad\text{for $s > 1.5$}\,.
\end{align*}
Since $\bbE''\cdot n = (b - \tb) + \teta''\cdot (n - \tn)$ and
$\teta\in L^\infty(0,T;H^{4.5}(\Gamma))$ by assumption, we find that
\begin{align}
& \int_0^t \|\bbE''\cdot n\|^2_{H^{2.5}(\Gamma)} ds \le C \int_0^t
\Big[\|b-\tb\|^2_{H^{2.5}(\Gamma)} +
\|n-\tn\|^2_{H^{2.5}(\Gamma)}\Big] ds \nonumber\\
\le&\ C t\int_0^t \Big[\|w_t\|^2_{H^1(\Omega)} +
\|w\|^2_{H^{2.5}(\Gamma)}\Big] ds + C \int_0^t
\Big[\|w_t\|^2_{L^2(\Omega)} + \|w\|^2_{H^{1.5}(\Gamma)} \Big] ds
\label{BH2p5estimate}\,.
\end{align}
%By (\ref{bbEH3p5estimate}) and $n\in L^\infty(0,T;
%H^{4.5}(\Gamma))$,
%\begin{align}
%& \int_0^t \Big[\|\bbE''\cdot n\|^2_{H^{2.5}(\Gamma)} + \|\bbE'''\cdot n\|^2_{H^{1.5}(\Gamma)}
%+ \|\bbE''''\cdot n\|^2_{H^{0.5}(\Gamma)}\Big] ds \label{BH2p5estimate}\\
%\le&\ C t\int_0^t \Big[\|w_t\|^2_{H^1(\Omega)} +
%\|w\|^2_{H^{2.5}(\Gamma)}\Big] ds + C \int_0^t
%\Big[\|w_t\|^2_{L^2(\Omega)} + \|w\|^2_{H^{1.5}(\Gamma)} \Big] ds
%\nonumber\,.
%\end{align}

\subsection{Estimates for $w_t$}
We study the time differentiated problem first. Time differentiating
(\ref{weq}) and then use $w_t$ as a test function, we find that
\begin{align}
&\ \ \frac{d}{dt} \Big[\frac{1}{2} \|w_t\|^2_{L^2(\Omega)} + 2
\||\eta_0'|^{-3/2}w'\cdot \eta'\|^2_{L^2(\Gamma)} +
\||\eta_0'|^{-3/2} w''\cdot n\|^2_{L^2(\Gamma)}
\Big] + \frac{\nu}{2} \|D_\eta w_t\|^2_{L^2(\Omega)} \nonumber\\
&= - \frac{\nu}{2} \langle (a_i^k)_t w_{,k}^j + (a_j^k)_t w^j_{,k},
(\Def_\eta w_t)_i^j \rangle - \frac{\nu}{2} \langle (\Def_\eta
w)_i^j, (a_i^k)_t w_{t,k}^i + (a_j^k)_t w_{t,k}^j\rangle \nonumber\\
&\ \ - \int_\Gamma |\eta_0'|^{-3}\Big[ (\tv'\cdot \bbE')
(\teta'\cdot
w')_t - (w'\cdot\eta') (v'\cdot w') - (\tv'\cdot \bbE') (v'\cdot w')\Big] dS \nonumber\\
&\ \ - \int_\Gamma |\eta_0'|^{-3}\Big[ - (w''\cdot n) (w''\cdot n_t)
+ [\tv''\cdot (n-\tn) + (\eta''\cdot n_t) -
(\teta''\cdot\tn_t)](w''_t\cdot n)\Big] dS \nonumber\\
&\ \ - \int_\Gamma |\eta_0'|^{-3}\Big[ 4 (g - \tg) (v'\cdot w'_t) + 2 (b - \tb) (w''_t\cdot n_t)\Big] dS \nonumber\\
&\ \ - \int_\Gamma |\eta_0'|^{-3}\Big[ 8 (\tv'\cdot\teta') (\bbE'\cdot w_t') + 2 (\tv''\cdot \tn + \teta''\cdot \tn_t) (w''_t \cdot (n-\tn))\Big] dS \nonumber\\
&\ \ - \int_\Gamma |\eta_0'|^{-3}\Big[ 4 (\tg - g_0) (w'\cdot w'_t) + 2 (\tb - b_0) (w''_t\cdot (n_t-\tn_t))\Big] dS \nonumber\\
&\ \ + \int_\Gamma |\eta_0'|^{-3} \Big[(g^{-1} - \tg^{-1})_t (b -
b_0) + (g^{-1} - \tg^{-1}) b_t\Big] g' (w_t'\cdot n) dS \label{L2H1wt}\\
&\ \ + \int_\Gamma |\eta_0'|^{-3} (g^{-1} - \tg^{-1}) (b - b_0)
\Big[g'_t (w_t'\cdot n) + g' (w_t'\cdot n_t)\Big] dS \nonumber\\
&\ \ + \int_\Gamma |\eta_0'|^{-3} \Big[(\tg^{-1})_t (b - \tb) +
\tg^{-1} (b - \tb)_t \Big] g' (w_t'\cdot n) dS \nonumber\\
&\ \ + \int_\Gamma |\eta_0'|^{-3} \tg^{-1} (b - \tb) \Big[g'_t
(w_t'\cdot n) + g'
(w_t'\cdot n_t)\Big] dS \nonumber\\
&\ \ + \int_\Gamma |\eta_0'|^{-3} \Big[(\tg^{-1})_t (\tb - b_0) +
\tg^{-1} \tb_t \Big] (g-\tg)' (w_t'\cdot n) dS \nonumber\\
&\ \ + \int_\Gamma |\eta_0'|^{-3} \tg^{-1} (\tb - b_0)
\Big[(g-\tg)'_t
(w_t'\cdot n) + (g-\tg)' (w_t'\cdot n_t)\Big] dS \nonumber\\
&\ \ + \int_\Gamma |\eta_0'|^{-3} \Big[(\tg^{-1})_t (\tb - b_0) + \tg^{-1} \tb_t\Big] \tg' (w_t'\cdot (n-\tn)) dS \nonumber\\
&\ \ + \int_\Gamma |\eta_0'|^{-3} \tg^{-1} (\tb - b_0) \Big[\tg'_t
(w_t'\cdot (n-\tn)) + \tg' (w_t'\cdot (n-\tn)_t)\Big] dS \nonumber\\
&\ \ - \int_\Gamma (\delta L_1)_t \cdot w_t dS + \langle (\delta
F)_t, w_t\rangle - \langle r_t, a_i^j w_{t,j}^i\rangle - \langle r,
(a_i^j)_t w_{t,j}^i \rangle \,.\nonumber
\end{align}
With $\tv_{tt} \in L^2(0,T;H^1(\Omega))$ (so that $v_{tt}$ has trace
on the boundary), similar to the computation of estimating $E_1$ in
page 48 of \cite{ChCoSh2006}, we find that
\begin{align*}
\Big|\langle r_t, a_i^j w_{t,j}^i\rangle + \langle r, (a_i^j)_t
w_{t,j}^i \rangle\Big| \le C_\epsilon \|r\|^2_{L^2(\Omega)} +
\epsilon \Big[\|r\|^2_{H^1(\Omega)} +
\|w_t\|^2_{H^1(\Omega)}\Big]\,.
\end{align*}
For $v$ and $\tv$ and the associated metric tensor, the second
fundamental form and the pressure in the space described in the
beginning of this section, by (\ref{deltaF}b),
\begin{align*}
|\langle (\delta F)_t, w_t\rangle| \le C_\epsilon \int_0^t
\|w\|^2_{{\mathcal V}^3(T)} ds + \epsilon \Big[\|w\|^2_{H^3(\Omega)}
+ \|w_t\|^2_{H^1(\Omega)}\Big]\,.
\end{align*}
By (\ref{deltaF}c) and interpolation inequalities,
\begin{align*}
\Big| \int_\Gamma (\delta L_1)_t \cdot w_t dS \Big| \le C t \int_0^t
\|w\|^2_{H^3(\Omega)} ds + C_\epsilon \|w_t\|^2_{L^2(\Omega)} +
\epsilon \|w_t\|^2_{H^1(\Omega)}\,.
\end{align*}
It is also clear that the first two terms (due to viscosity) on the
right-hand side is bounded by $C \|w\|_{H^1(\Omega)}
\|w_t\|_{H^1(\Omega)}$, and by Young's inequality,
\begin{align*}
& \Big|\langle (a_i^k)_t w_{,k}^j + (a_j^k)_t w^j_{,k}, a_i^k
w_{t,k}^i + a_j^k w_{t,k}^j\rangle + \langle a_i^k w_{,k}^j + a_j^k
w^j_{,k}, (a_i^k)_t w_{t,k}^i + (a_j^k)_t w_{t,k}^j\rangle\Big| \\
\le&\ C_\epsilon \int_0^t \|w_t\|^2_{H^1(\Omega)} ds + \epsilon
\|w_t\|^2_{H^1(\Omega)} \,.
\end{align*}
For those terms having $w_t'$ in the integrands, %since
%\begin{align*}
%& \|\eta - \teta\|^2_{H^{1.5}(\Gamma)} + \|w\|^2_{H^{1.5}(\Gamma)} +
%\|g - \tg\|^2_{H^{1.5}(\Gamma)} + \|\tb - b_0\|^2_{H^{1.25}(\Gamma)}\|(g - \tg)_t\|^2_{H^{1.5}(\Gamma)} \\
%+&\ \|b-\tb\|^2_{H^{0.5}(\Gamma)} +
%\underline{\|(b-\tb)_t\|^2_{H^{0.5}(\Gamma)}} (???) +
%\|n-\tn\|^2_{H^{0.5}(\Gamma)} +
%\|(n-\tn)_t\|^2_{H^{0.5}(\Gamma)} \\
%\le&\ C_\epsilon \int_0^t \|w\|^2_{{\mathcal V}^3(T)} ds + (\epsilon + C
%t^\delta) \|w\|^2_{H^3(\Omega)} \,,
%\end{align*}
%using $H^{0.5}(\Gamma)$-$H^{-0.5}(\Gamma)$ duality pairing
following the same procedure of estimating $L_4$, we find that those
terms are bounded by
\begin{align*}
C_\epsilon \Big[\int_0^t \|w\|^2_{{\mathcal V}^3(T)} ds +
\|w''\|^2_{H^1(\Omega_1')}\Big] + \epsilon
\|w_t\|^2_{H^1(\Omega)}\,.
\end{align*}
The terms having $w_t''$ inside the integrands are
\begin{align*}
M_1 &= \int_\Gamma |\eta_0'|^{-3} (b-\tb) (w''_t\cdot n_t)\Big] dS\,,\\
M_2 &= \int_\Gamma |\eta_0'|^{-3} \Big[(\eta''\cdot n_t) -
(\teta''\cdot\tn_t)\Big](w''_t\cdot n) dS\,,\\
M_3 &= \int_\Gamma \Big[F_1(\eta_0,\eta,\teta) (n-\tn) w_t'' +
F_2(\eta_0,\eta,\teta) (\tb - b_0) (n-\tn)_t w_t''\Big] dS\,.
\end{align*}
Following the same procedure of estimating $L_3$, we find that
\begin{align*}
|M_1| + |M_3| \le&\ C_\epsilon \Big[\int_0^t \|w\|^2_{H^3(\Omega)}
ds + t^\delta \|w\|^2_{H^3(\Omega)} + \|r\|^2_{H^1(\Omega)} + \|w\|^2_{H^2(\Omega)}\Big] \\
& + \epsilon \|w_t\|^2_{H^1(\Omega)}
\end{align*}
and hence
\begin{align*}
\Big|\int_0^t (M_1 + M_3) ds\Big| \le&\ C_\epsilon \Big[ t^\delta
\int_0^t \|w\|^2_{H^3(\Omega)} ds + \int_0^t \|w_t\|^2_{L^2(\Omega)}
ds \Big]\,.
\end{align*}
For $M_2$, by interpolations and (\ref{bbEH3p5estimate}), its time
integral satisfies
\begin{align*}
\Big|\int_0^t M_2 ds\Big| \le&\ C_\epsilon \int_0^t
\Big[\|\bbE\|^2_{H^{3.25}(\Gamma)} + \|w''\|^2_{H^1(\Omega_1')}\Big]
ds +
\epsilon \int_0^t \|w_t\|^2_{H^1(\Omega)} ds \\
\le&\ C_\epsilon \int_0^t \|\bbE\|^2_{H^{2.5}(\Gamma)} ds + C_\epsilon \int_0^t \|w''\|^2_{H^1(\Omega_1')} ds\\
& + \epsilon \int_0^t \|\bbE\|^2_{H^{3.5}(\Gamma)} ds + \epsilon \int_0^t \|w_t\|^2_{H^1(\Omega)} ds \\
\le&\ C_\epsilon \int_0^t \|w''\|^2_{H^1(\Omega_1')} ds + \epsilon
\int_0^t \|w_t\|^2_{H^1(\Omega)} ds\,.
\end{align*}
All the remaining terms %have either $w''$ in pair or $\bbE'$
%together with $w'$, so they
can be bounded by
\begin{align*}
C \Big[\|w\|^2_{H^2(\Gamma)} + \int_0^t \|w\|^2_{H^(\Gamma)}
ds\Big]\,.
\end{align*}
Time integrating (\ref{L2H1wt}), choosing $\epsilon>0$ and $T>0$
small enough, by (\ref{wrregularity1}) and (\ref{wrregularity3}) we
find that
\begin{align}
Y(T) + \int_0^T \|w_t\|^2_{H^1(\Omega)} ds %\le&\ C \int_0^T
%\|w_t\|^2_{L^2(\Omega)} ds + C \int_0^T \|w''\|^2_{H^1(\Omega_1')}
%ds \nonumber\\
\le&\ C \int_0^T Y(T) ds + C \int_0^T \|w''\|^2_{H^1(\Omega_1')} ds
\label{Yineq}
\end{align}
where
\begin{align*}
Y(T) = \sup_{0\le t\le T} \Big[\|w_t\|^2_{L^2(\Omega)} +
\|w'\cdot\eta'\|^2_{L^2(\Gamma)} + \|w''\cdot
n\|^2_{L^2(\Gamma)}\Big]\,.
\end{align*}

\subsection{Estimates for $w''$}
\noindent Let $\phi = (\zeta_1^2 w'')''$ in (\ref{weq}), then
\begin{align}
&\ \ \frac{d}{dt} \Big[\frac{1}{2}\|\zeta_1 w_t''\|^2_{L^2(\Omega)}
+ 2\||\eta_0'|^{-3/2}(\bbE' \cdot \eta')''\|^2_{L^2(\Gamma)} +
2\||\eta_0'|^{-3/2}(\bbE' \cdot \teta')''\|^2_{L^2(\Gamma)} \nonumber\\
&\qquad + \||\eta_0'|^{-3/2}(\bbE''\cdot n)''\|^2_{L^2(\Gamma)}\Big]
+ \frac{\nu}{2} \|\zeta_1 D_\eta w''\|^2_{L^2(\Omega)} \nonumber\\
&= - \nu\int_{\Omega} \zeta_1 \Big[2 (a_i^{k\prime} w_{,k}^{j\prime}
+ a_j^{k\prime} w_{,k}^{i\prime}) + (a_i^{k\prime\prime} w_{,k}^j +
a_j^{k\prime\prime} w^i_{,k})\Big] a_i^k \zeta_1 w_{,k}^{j\prime\prime} dx \nonumber\\
&\ \ - 2 \nu \int_{\Omega} \zeta_1 \Big[(a_i^{k\prime} w^j_{,k} +
a_j^{k\prime} w^i_{,k}) + (a_i^k w_{,k}^{j\prime} + a_j^k
w_{,k}^{i\prime})\Big] a_i^{k\prime} \zeta_1
w_{,k}^{j\prime\prime} dx \nonumber\\
&\ \ + 4 \int_\Gamma |\eta_0'|^{-3} (g - \tg) \Big[ (\eta^{(5)}
\cdot\bbE')_t + 4 (\eta^{(4)}\cdot\bbE'')_t + 6
(\eta'''\cdot\bbE''')_t + 4 (\eta''\cdot\bbE^{(4)})_t\Big] dS
\nonumber\\
&\ \ - 4 \int_\Gamma \Big[(|\eta_0'|^{-3})' (\bbE'\cdot \eta')' + (|\eta_0'|^{-3})'' (\bbE'\cdot\eta')\Big] (\eta'\cdot \bbE')_t'' dS \nonumber\\
&\ \ - 4 \int_\Gamma \Big[(|\eta_0'|^{-3})' (\bbE'\cdot \teta')' + (|\eta_0'|^{-3})'' (\bbE'\cdot\teta')\Big] (\teta'\cdot \bbE')_t'' dS \nonumber\\
&\ \ - 4 \int_\Gamma |\eta_0'|^{-3} (\bbE'\cdot\teta')'' (\bbE'\cdot
\bbE')_t'' dS
- \int_\Gamma |\eta_0'|^{-3} \Big[\teta''\cdot (n-\tn)\Big]'' (n\cdot \bbE'')''_t dS \label{L2H3w}\\
&\ \ - \int_\Gamma \Big[2(b-\tb)' (|\eta_0'|^{-3} n)' + (b - \tb)(|\eta_0'|^{-3} n)'' \Big]\cdot w^{(4)} dS \nonumber\\
&\ \ + \int_\Gamma |\eta_0'|^{-3} (b - \tb)'' \Big[(\bbE''\cdot n_t)'' + 2 w'''\cdot n' + w''\cdot n''\Big] dS \nonumber\\
&\ \ - 4 \int_\Gamma \Big[|\eta_0'|^{-3}(\tg - g_0)
\bbE'\Big]''\cdot w''' dS - 2 \int_\Gamma \Big[|\eta_0'|^{-3} (\tb -
b_0) (n-\tn)\Big]'' \cdot w^{(4)} dS \nonumber\\
&\ \ - \int_\Gamma \Big[|\eta_0'|^{-3} \Big(\bbE'\cdot \eta' +
\teta'\cdot\bbE'\Big) g^{-1} \tg^{-1} (b-b_0) g'
n\Big]''\cdot w''' dS \nonumber\\
&\ \ + \int_\Gamma \Big[|\eta_0'|^{-3} \tg^{-1}\Big(\bbE''\cdot
n + \teta''\cdot (n-\tn)\Big) g' n\Big]''\cdot w''' dS \nonumber\\
&\ \ + \int_\Gamma \Big[|\eta_0'|^{-3} \tg^{-1} (\tb - b_0)
\Big(\bbE'\cdot \eta' + \teta'\cdot\bbE'\Big)' n \Big]''
\cdot w''' dS \nonumber\\
&\ \ + \int_\Gamma \Big[|\eta_0'|^{-3} \tg^{-1} (\tb - b_0) \tg' (n
- \tn) \Big]''\cdot w''' dS \nonumber\\
&\ \ - \int_\Gamma (\delta L_1)\cdot w'''' dS - \langle r, a_i^j
(\zeta_1^2 w^{i\prime\prime})''_{,j} \rangle + \langle \delta F,
(\zeta_1^2 w'')'' \rangle \nonumber\,.
\end{align}
As the estimate of $K_1$ in Section \ref{Kestimate}, the first two
integrals (due to viscosity) on the right-hand side can be bounded
by $C \|w\|_{H^{2.5}(\Omega)} \|w\|_{H^3(\Omega)}$, and by Young's
inequality,
\begin{align*}
&\ \ \Big|\int_{\Omega} \zeta_1 \Big[2 (a_i^{k\prime}
w_{,k}^{j\prime} + a_j^{k\prime} w_{,k}^{i\prime}) +
(a_i^{k\prime\prime} w_{,k}^j + a_j^{k\prime\prime} w^i_{,k})\Big]
a_i^k \zeta_1
w_{,k}^{j\prime\prime} dx\Big| \\
&\ + \Big|\int_{\Omega} \zeta_1 \Big[(a_i^{k\prime} w^j_{,k} +
a_j^{k\prime} w^i_{,k}) + (a_i^k w_{,k}^{j\prime} + a_j^k
w_{,k}^{i\prime})\Big] a_i^{k\prime} \zeta_1 w_{,k}^{j\prime\prime}
dx \Big| \\
&\le C_\epsilon \|w\|^2_{H^2(\Omega)} + \epsilon
\|w\|^2_{H^3(\Omega)} \,.
\end{align*}
Similar to the computation of estimating $D_1$ in page 48 of
\cite{ChCoSh2006}, we find that
\begin{align*}
|\langle r, a_i^j (\zeta_1^2 w^{i\prime\prime})''_{,j} \rangle| \le
C_\epsilon t \int_0^t \|w\|^2_{H^3(\Omega)} + \epsilon
\|r\|^2_{H^2(\Omega)} \,.
\end{align*}
By (\ref{deltaF}a), (\ref{deltaF}c) and interpolation inequalities,
\begin{align*}
|\langle \delta F, (\zeta_1^2 w'')'' \rangle| + \Big|\int_\Gamma
(\delta L_1)\cdot w'''' dS\Big| \le C_\epsilon t\int_0^t
\|w\|^2_{H^3(\Omega)} ds + \epsilon \|w\|^2_{H^3(\Omega)}\,.
\end{align*}
By $H^{1.5}(\Gamma)$-$H^{-1.5}(\Gamma)$ or
$H^{0.5}(\Gamma)$-$H^{-0.5}(\Gamma)$ duality pairing,
\begin{align*}
& \Big| \int_\Gamma |\eta_0'|^{-3} (g - \tg) \Big[ (\eta^{(5)}
\cdot\bbE')_t + 4 (\eta^{(4)}\cdot\bbE'')_t + 6
(\eta'''\cdot\bbE''')_t + 4 (\eta''\cdot\bbE^{(4)})_t\Big] dS \Big| \\
\le&\ C \|g - \tg\|_{H^{1.5}(\Gamma)} \Big[ \|w\|_{H^{2.5}(\Gamma)}
+ \|\bbE\|_{H^{2.5}(\Gamma)} \Big]
\end{align*}
and
\begin{align*}
&\ \ \ \Big|\int_\Gamma \Big[(|\eta_0'|^{-3})' (\bbE'\cdot \eta')' +
(|\eta_0'|^{-3})'' (\bbE'\cdot\eta')\Big] (\eta'\cdot \bbE')_t'' dS
\Big| \\
&\ + \Big|\int_\Gamma \Big[(|\eta_0'|^{-3})' (\bbE'\cdot \eta')' +
(|\eta_0'|^{-3})'' (\bbE'\cdot\eta')\Big] (\eta'\cdot \bbE')_t'' dS
\Big| \\
&\le C \Big[\|\bbE'\cdot\eta'\|_{H^{1.5}(\Gamma)} +
\|\bbE'\cdot\teta'\|_{H^{1.5}(\Gamma)}\Big]\Big[\|\bbE\|_{H^{2.5}(\Gamma)}
+ \|w\|_{H^{2.5}(\Gamma)}\Big]\,.
\end{align*}
Standard H$\ddot{\text{o}}$older's inequality shows that
\begin{align*}
\Big|\int_\Gamma |\eta_0'|^{-3} (\bbE'\cdot\teta')'' (\bbE'\cdot
\bbE')_t'' dS\Big| \le C \Big[\|\bbE'\cdot \teta'\|^2_{H^2(\Gamma)}
+ \|\bbE\|^2_{H^3(\Gamma)} \Big]\,.
\end{align*}
Similar to the estimates of $K_3$, $K_4$ and $K_5$, by
(\ref{bgsmall}), (\ref{bL2H25estimate}) and (\ref{BH2p5estimate}),
we find that the time integral of the remaining terms can be bounded
by
\begin{align*}
(C_\epsilon t + \epsilon) \int_0^t \Big[\|w_t\|^2_{H^1(\Omega)} +
\|w\|^2_{H^{2.5}(\Gamma)}\Big] ds + C_\epsilon \int_0^t
\Big[\|w_t\|^2_{L^2(\Omega)} + \|w\|^2_{H^{1.5}(\Gamma)} \Big] ds\,.
\end{align*}
Time integrating (\ref{L2H3w}), we find that
\begin{align}
Z(T) + \int_0^T \|w''\|^2_{H^1(\Omega_1')} ds \le&\ (C_\epsilon T +
\epsilon) \int_0^T \Big[\|w_t\|^2_{H^1(\Omega)} +
\|w\|^2_{H^{2.5}(\Gamma)}\Big] ds \label{Zineq}\\
& + C_\epsilon \int_0^T (Y(t) + Z(t)) dt\,, \nonumber
\end{align}
where
\begin{align*}
Z(T) =& \sup_{0\le t\le T} \Big[\|w''\|^2_{L^2(\Omega_1')} +
\|\bbE'\cdot\eta'\|^2_{H^2(\Gamma)} +
\|\bbE'\cdot\teta'\|^2_{H^2(\Gamma)} + \|\bbE''\cdot
n\|^2_{H^2(\Gamma)} \Big] \,.
\end{align*}
Combining (\ref{Yineq}) and (\ref{Zineq}), choosing $\epsilon>0$
small enough and then $T>0$ small enough, we find that
\begin{align*}
Y(T) + Z(T) + \int_0^T \Big[\|w_t\|^2_{H^1(\Omega)} +
\|w''\|^2_{H^1(\Omega_1')} \Big] ds \le C \int_0^T (Y(t) + Z(t))
dt\,.
\end{align*}
By the Gronwall inequality, $Y$ and $Z$ are identical to zero, which
shows $v = \tv$, and hence the solution $v\in {\mathcal V}^5(T)$ to
(\ref{solidshelleq}) is unique.


\begin{thebibliography}{50}

\bibitem{AuBe} F.~Auricchio, L. Beir${\tilde{\rm{a}}}$o da Veiga
and C. Lovadina, {\scshape Remarks on the asymptotic behaviour of
Koiter shells,} Computers and Structures, {\bf 80} (2002), 735-745

\bibitem{DaVe} H. Beir\~{a}o. da Veiga, {\scshape On the existence of strong
solutions to a coupled fluid-structure evolution problem,} J. Math.
Fluid Mech.,  {\bf 6}  2004,  21--52.

\bibitem{Bo2005} M.~Boulakia, {\scshape Existence of weak solutions for an interaction
problem between an elastic structure and a compressible viscous
fluid,} J. Math. Pures Appl. (9) {\bf 84} (2005),  no. 11,
1515--1554.

\bibitem{ChDeEsGr} A.~Chambolle, B.~Desjardins, M.J.~Esteban, C.~Grandmont,
{\scshape Existence of weak solutions for an unsteady fluid-plate
interaction problem,} Preprint.

\bibitem{ChCoSh2006} C.H.~Cheng, D.~Coutand, S. Shkoller, {\scshape Navier-Stokes equations interacting with a nonlinear elastic
shell,} Preprint.

\bibitem{PGC} P.G.~Ciarlet, {\scshape An introduction to
Differential Georetry with Applications to Elasticity}, Springer
2005

\bibitem{PGC1} P.G.~Ciarlet, {\scshape Introduction to linear shell theory,} Series in Applied
Mathematics (Paris), vol. 1, Gauthier-Villars, Editions
Scientifiques et M edicales Elsevier, Paris, 1998.

\bibitem{CoSh2002} D.~Coutand and S. Shkoller, {\scshape Unique solvability of the free-boundary
Navier-Stokes equations with surface tension,}

\bibitem{CoSh2005} D.~Coutand and S. Shkoller, {\scshape On the motion of an elastic solid
inside of an incompressible viscous fluid,} to appear in Arch.
Rational Mech. Anal.

\bibitem{CoSh2006} D.~Coutand and S. Shkoller, {\scshape On the interaction between quasilinear
elastodynamics and the Navier-Stokes equations,}

\bibitem{DeEs} B.~Desjardins, M.J.~Esteban,
 {\scshape  Existence of weak solutions for the motion of rigid bodies in a viscous fluid,} Arch. Rational Mech. Anal.,
{\bf 146} (1999), 59--71.

\bibitem{DeEsGrTa} B.~Desjardins, M.J.~Esteban, C.~Grandmont, P.~Le Tallec,
 {\scshape Weak solutions for a fluid-structure interaction problem,} Rev. Mat.
Complut., {\bf 14} (2001), 523--538.

\bibitem{LCE1} L.C.~Evans, {\scshape Partial Differential Equations},
Graduate Studies in Mathematics, {\bf 19} American Mathematical
Society, Providence, RI, 1998.

\bibitem{FlOr} F.~Flori, P.~Orenga, {\scshape Fluid-structure interaction: analysis of a 3-D compressible model,} Ann. Inst. H. Poincar\'e Anal. Non Lin\'eaire, {\bf 17} (2000), 753-777.

\bibitem{GPG1} G.P.~Galdi, {\scshape An Introduction to the Mathematical Theory of the
Navier-Stokes Equations Volume I,} Springer Tracts in Natural
Philosophy, Vol {\bf 38}.

\bibitem{GeKrMa1996} Z.~Ge, H.P. Kruse and J.E. Marsden, {\scshape The limits of Hamiltonian
structures in three-dimensional elasticity, shells, and rods,} J.
Nonlinear Sci. Vol. {\bf 6} (1996), 19-57.

\bibitem{EG} E.~Givelberg, {\scshape Modeling Elastic Shells Immersed in Fluid,} Comm. Pure Appl.
Math., {\bf 57} (2004), no. 3, 283-309.

\bibitem{GrMa} C.~Grandmont, Y.~Maday,
 {\scshape  Existence for unsteady fluid-structure interaction problem,} Math. Model. Numer. Anal.,
{\bf 34} (2000), 609--636.

\bibitem{LePeLa} R.J.~Leveque, C.S. Peskin and P.D. Lax, {\scshape Solution of
a two-dimensional cochlea model with fluid viscosity,} SIAM J. Appl.
Math., {\bf 45} (1985), no. 3, 450-464.

\bibitem{LiWa} C. Liu, N.J. Walkington, {\scshape An Eulerian description of fluids containing visco-elastic particles}, Arch. Rational Mech. Anal., {\bf 159} (2001), 229-252.

\bibitem{Serre} D.~Serre,
 {\scshape  Chute libre d'un solide dans un fluide visqueux incompressible: Existence,} Japan J. Appl. Math.,
{\bf 4} (1987), 33--73.

\bibitem{Wein} H.F.~Weinberger, {\scshape Variational properties of steady fall in Stokes flow,} J. Fluid Mech., {\bf 52} 1972, 321--344.

\end{thebibliography}
\end{document}